\documentclass[pdflatex,
iicol,
sn-mathphys-num]{sn-jnl}

\usepackage{graphicx}%
\usepackage{multirow}%
\usepackage{amsmath,amssymb,amsfonts}%
\usepackage{amsthm}%
\usepackage{mathrsfs}%
\usepackage[title]{appendix}%
\usepackage{xcolor}%
\usepackage{textcomp}%
\usepackage{manyfoot}%
\usepackage{booktabs}%
\usepackage{algorithm}%
\usepackage{algorithmicx}%
\usepackage{algpseudocode}%
\usepackage{listings}%
\usepackage{tabularx}
\newcolumntype{Y}{>{\centering\arraybackslash}X}
\usepackage{mathtools} 
\usepackage{microtype}
\usepackage{bm}  
\usepackage{stmaryrd} 
\usepackage{subcaption}
\usepackage{siunitx}
\usepackage{tikz}
\usetikzlibrary{calc}
\usetikzlibrary{patterns}
\usetikzlibrary{arrows.meta, shapes.misc}

\tikzset{
  cross/.style={draw, label={[font=\large]center:$\times$}, circle,
    inner sep=1mm, outer sep=0pt},
  dot/.style={draw, circle, inner sep=1mm, outer sep=0, label={[fill=black, circle, inner sep=.5mm]center:}},
  point/.style={fill=black, circle, inner sep=1pt, node contents={}},
}

\graphicspath{ {./figures/} }

\definecolor{plotcolor1}{HTML}{1F46A1}
\definecolor{plotcolor2}{HTML}{6E3F36}

\theoremstyle{thmstyleone}%
\theoremstyle{thmstyletwo}%

\theoremstyle{thmstylethree}%

\newcommand{\ltwo}[1]{\textrm{L}^2 (#1)}
\newcommand{\hcurl}[1]{\textrm{H} (\textrm{curl};#1)}
\newcommand{\hzcurl}[1]{\textrm{H}_0(\textrm{curl};#1)}
\newcommand{\hone}[1]{\textrm{H}^1(#1)}
\newcommand{\hzone}[1]{\textrm{H}^1_0(#1)}

\DeclareMathOperator{\divop}{\mathrm{div}}
\DeclareMathOperator{\curlop}{\nabla \times}
\DeclareMathOperator{\gradop}{\nabla}
\DeclareMathOperator{\curlopt}{\nabla_t \times}
\DeclareMathOperator{\gradopt}{\nabla_t}

\begin{document}

\title[wpbc]{A Waveguide Port Boundary Condition Based on Approximation Space Restriction for Finite Element Analysis}


\author*[1]{\fnm{Francisco} \sur{T. Orlandini}}\email{f102348@dac.unicamp.br}

\author[2]{\fnm{Philippe} \sur{R. B. Devloo}}\email{phil@unicamp.br}

\author[1]{\fnm{Hugo} \sur{E. Hernández-Figueroa}}\email{hugo@unicamp.br}

\affil*[1]{\orgdiv{DECOM}, \orgname{Faculdade de Engenharia Elétrica e Computação}, \orgaddress{Universidade Estadual de Campinas, Brazil}}

\affil[2]{\orgdiv{Departamento de Estruturas}, \orgname{Faculdade de Engenharia Civil, Arquitetura e Unicamp}, \orgaddress{Universidade Estadual de Campinas, Brazil}}


\abstract{A Waveguide Port Boundary Condition (WPBC)
    based on the restriction of the approximation space is presented
    in the context of Finite Element Analysis.
    As well as reducing the computational domain in the same manner as the traditional
    WPBC,
    the proposed scheme further reduces the degrees of freedom at the waveguide ports,
    simplifies the implementation and seamlessly provides post-processing results such
    as the reflected power in each waveguide mode.
    The boundary condition is thoroughly derived, and numerical examples are used
    as a support for the discussion on topics such as the needed number of modes
    to be employed at a waveguide port. Finally, a nanograting-based plasmonic sensor
    is analysed to illustrate further possibilities of the scheme.}

\keywords{FEM, optical waveguides, reflectionless boundary conditions}



\maketitle

\section{Introduction}\label{sec:intro}
From submarine communications cables to Radio Frequency Identification systems,
electromagnetics waveguides are of great importance in several fields of electrical engineering.
Moreover, advancements in \emph{Photonic Crystal Fibres} (PCFs) have allowed the guiding and confining light in
dimensions previously unattainable until recent decades.
Therefore, the class of \emph{waveguide-fed} optical devices is an ever growing class, fuelling thus the
need for an accurate numerical analysis of complex waveguide structures.

One class of numerical method that has been employed for this end is the \emph{Finite Element Method} (FEM),
which, among its desirable features, present a versatility to deal with complex geometries.
The frequency domain analysis of these devices is a task often needed in the design process, as one could
then compute reflectivity/transmission spectra, field distributions and signals coupling.

However, FEM schemes require a \emph{closed} computational domain, and the truncation of open domains has been a subject
of study for several years, as show in \citet{merewether1971}, in the context of Finite Difference schemes.
In the context of waveguide-fed structures, one special case of open domains of particular interest are the
\emph{waveguide ports}, in which one can inject a given incident wave and absorb the reflected
\emph{waveguide modes} from the structure, and it its thus important to truncate the domain appropriately
at the ports.

The \emph{Perfectly Matching Layers} (PMLs), introduced in \citet{berenger1994}, are a popular choice
among the Computational Electromagnetics community. The PMLs can be interpreted as a fictitious region
around the computational domain of interest that are able to attenuate incident waves, and are often placed
behind a waveguide port. However, not only they result in an expansion of the domain (and, naturally,
in an increase of degrees of freedom in a region away from the domain of interest) but adjusting the PMLs
for absorbing evanescent waves as well as propagating waves can be a challenging task\citep{ouldagha2008},
specially over a range of frequencies.

Among other approaches\citep{davanco2001,rubiomercedes2004},
of special significance for the present work is the
\emph{Waveguide Port Boundary Condition} (WPBC) \citep{jin1994}, which
is based on an orthogonality relationship between waveguide modes and can accurately absorb the reflected
modes. In \citet{jin1994} the WPBC is described for rectangular waveguides, and in \citet{jin2014} a more
general form is presented.

However, there are a couple of issues, the first being the orthogonality relationship
between the waveguide modes. Taking a step-index optical fibre as an example of
a waveguide with an open cross section, one can verify that only the guided modes
present an exponential decay away from the core. With radiation modes, which are
also needed to represent a reflection, this does not happen. How could one
derive an orthogonality relationship for a mode that is unbounded in the spatial domain
and thus is not even square-integrable?

This property of some of the physical solutions translate to the computational domain
as well, given that a lossless open waveguide cross section truncated with PMLs
result in a non-self adjoint problem, and the orthogonality must be carefully studied,
as recently discussed in \citet{sauvan2022}.

At last, even if the orthogonality has been dealt with,
the traditional form of the WPBC involves double integrals of the boundary
term, which results in a more complicated implementation.

A novel form of the WPBC is presented in this paper, wherein the \emph{dofs}
at the boundary are associated with the waveguide modes due to an approximation space restriction.
Its expressions are derived and the orthogonality relationship of the eigenmodes
is discussed for lossy, transverse-anisotropic waveguides,
allowing for its application in a general waveguide. The paper is organised as follows:
Section \ref{sec:pre} presents and discusses the weak formulations and the necessary Hilbert spaces.
Section \ref{sec:wpbc} dwells on the orthogonality relationship between the FEM solutions
of the waveguide analysis and from thereon develops the
new WPBC formulation with approximation space restriction, along with a discussion on
aspects related to its implementation. Finally, in Section \ref{sec:res}, numerical
results are presented and strategies for the application of the WPBC are discussed,
including a discussion on a strategy to determine the needed number of modes
for a reflectionless boundary condition in a more complex scenario.

\section{Preliminaries}\label{sec:pre}
\subsection{Hilbert spaces}\label{sec:hilb}
The relevant function spaces for the FEM formulations are the following:
\noindent\begin{align}
  &\begin{multlined}\label{eq:fem-l2-space-def}
    \ltwo{\Omega} =\left\lbrace\, \phi: \Omega \rightarrow \mathbb{C}\,\middle| \,\int_\Omega \phi\cdot \phi^*\,d\Omega < \infty \,\right\rbrace\text{,}
  \end{multlined}\\
  &\begin{multlined}\label{eq:fem-h1-space-def}
    \hone{\Omega} =\left\lbrace\, \phi \in \ltwo{\Omega}\,\middle| \, \gradop \phi \in [\ltwo{\Omega}]^d\,\right\rbrace\text{,}
  \end{multlined}\\
  &\begin{multlined}\label{eq:fem-hcurl-space-def}
    \hcurl{\Omega} = \\
    \left\lbrace\, \bm{\phi} \in [\ltwo{\Omega}]^d\,\middle| \, \curlop \bm{\phi} \in [\ltwo{\Omega}]^{d'}\,\right\rbrace\text{.}
    \end{multlined}
\end{align}

The $\ltwo{\Omega}$ space refers to the scalar functions which are square-integrable in the domain $\Omega$,
$\hone{\Omega}$ is a subspace of $\ltwo{\Omega}$ in which the gradient is also square-integrable and, finally,
$\hcurl{\Omega}$ denotes the subspace of $[\ltwo{\Omega}]^d$, where $d$ is the dimension of $\Omega$,
for which the curl is well defined and it lies in $[\ltwo{\Omega}]^{d'}$, with $d'=2d-3$ for $d=2,3$.

For $\hone{\Omega}$ and $\hcurl{\Omega}$ we also have the associated \emph{trace} spaces.
Restricting the discussion to quasi-polynomial spaces for clarity,
the trace of a function with support on $\Omega$ is associated with its restriction on the boundary $\Gamma$.
For $\hone{\Omega}$, the trace is merely the restriction of the function on $\Gamma$,
while for $\hcurl{\Omega}$ the trace is the \emph{tangential component} of the function. Formally, one can write
\begin{align}
  &\begin{multlined}
    \mathrm{tr}: \phi \in \hone{\Omega} \rightarrow \mathrm{tr}\,\phi =\phi\vert_{\Gamma}\in \ltwo{\Gamma}
  \end{multlined}\\
  &\begin{multlined}
    \mathrm{tr}: \bm{\phi} \in \hcurl{\Omega} \rightarrow\\
    \mathrm{tr}\,\bm{\phi} =
    \bm{n}\times\bm{\phi}\vert_{\Gamma}\in \ltwo{\Gamma}^{d'}\text{,}
  \end{multlined}
\end{align}
and, having the traces defined
\footnote{Actually, the traces reside in subspaces of the corresponding $\ltwo{\Omega}$ spaces,
  but we avoid defining fractional Hilbert spaces given the scope of the paper.},
present the subspaces
\begin{align}
  &\begin{multlined}\label{eq:fem-hz1-space-def}  
    \hzone{\Omega} =
    \left\lbrace\, \phi \in \hone{\Omega}\,\middle| \,
    \mathrm{tr}\,\phi\vert_{\Gamma_d}=0\right\rbrace\text{,}
\end{multlined}\\
&\begin{multlined}\label{eq:fem-hzcurl-space-def}
  \hzcurl{\Omega} =\\
  \left\lbrace\, \bm{\phi} \in \hcurl{\Omega}^d\,\middle| \,
    \mathrm{tr}\,\bm{\phi}\vert_{\Gamma_d} = 0\right\rbrace\text{,}
  \end{multlined}
\end{align}
which are used throughout this work for imposing homogeneous Dirichlet boundary conditions on the
$\Gamma_d\subseteq \Gamma$ boundaries.

\subsection{Weak formulations}\label{sec:weak}
All the weak formulations presented in this paper are derived from the inhomogeneous Helmholtz
wave equation, shown in Equation \ref{eq:helmholtz}.
The derivations are presented here for the complex electric field, however, they
can also be formulated using the magnetic field due to the duality of the Maxwell's equations.

\begin{equation}\label{eq:helmholtz}
  \curlop\left[\bm{\mu}^{-1}_r\left(\curlop\bm{E}\right)\right]-k_0^2\bm{\epsilon}_ r\bm{E}
  = \bm{F}\text{.}
\end{equation}

The relative electric permittivity and magnetic permeability are shown here in the tensor form
as to contemplate anisotropy, and are denoted, respectively, by $\bm{\epsilon}_r$ and $\bm{\mu}_r$.
$\bm{F}$ is a general source term and $k_0$ represents the wavenumber in vacuum.

The framework we are interested in is to perform the modal analysis of a waveguide which is embedded in a
domain in which we are interested to analyse the wave propagation, or \emph{scattering}, resulting from the
incidence of one or more waveguide modes. This section is organised in the same fashion, discussing first the modal
analysis weak formulation followed by the scattering formulation.

\subsubsection{Modal analysis formulation}\label{sec:modal}

In the context of modal analysis, a waveguide is a structure invariant in the $z$-direction
\footnote{The scope of this work does not contemplate waveguides with periodicity in the propagation direction.},
which is said to be the \emph{propagation} direction.
Moreover, the constitutive parameters $\bm{\epsilon}_r$ and $\bm{\mu}_r$ are now restricted to transverse anisotropy
and the source term $\bm{F}$ is zero.

Assuming that the waveguide has a (possibly unbounded) cross-section $\Omega$ contained in
the $xy$-plane, the modal analysis problem is posed as finding electric fields of the form

\begin{equation}
  \label{eq:e-modal-analysis}
  \bm{E}(x,y,z) = \bm{E}(x,y)e^{-j\beta z}
\end{equation}
satisfying Equation \ref{eq:helmholtz} with $\bm{F}=0$. $\beta$ is the \emph{propagation constant} of the
waveguide mode.

Given the invariance in the propagation direction and the known $z$ dependency of the electric field, the
modal analysis problem can be solved only in the two-dimensional cross section $\Omega$. The weak form is
essentially the one presented in \citet{lee1991}, but supporting also transverse anisotropy.

The weak formulation is obtained by applying Equation \ref{eq:e-modal-analysis} in \ref{eq:helmholtz}
and splitting the electric field as $\bm{E}=\bm{E}_t+E_z\hat{\bm{z}}$ and the differential operator as
$\nabla = \nabla_t+ \frac{\partial}{\partial z} \hat{\bm{z}}$.
After some algebra, one can use the transformation of variables \citep{lee1991}:

\begin{equation}\label{eq:fem-hcurl-formulation-transf}
  \bm{e}_t=\beta\bm{E}_t\qquad e_z=-jE_z\text{,}
\end{equation}
and obtain the following weak formulation:

Find non-trivial
$\left(\beta^2,\bm{e}_t, e_z\right)$ $\in $ $(\mathbb{C} \times \hzcurl{\Omega} \times \hzone{\Omega})$
such that:
\begin{multline}\label{eq:wgma-var-final}
  \int_\Omega
  \left[\bm{\mu}_{r}^{-1} \left( \curlopt \bm{e}_t\right)\right]\cdot \left( \curlopt \bm{\phi} \right)^*
  - k_0^2 \bm{\epsilon}_r\bm{e}_t\cdot\bm{\phi}^*\mathrm{d}\Omega\\*
  +\beta^2\int_\Omega
  \left\{\bm{\mu}_{r}^{-1} \left[\left( \nabla_t e_z +\bm{e}_t \right)\times\hat{\bm{z}}\right]\right\}
 \cdot
  \left[\left( \nabla_t \phi +\bm{\phi} \right)\times\hat{\bm{z}}\right]^*\\*
  -k_0^2\left[\bm{\epsilon}_{r} \left(e_z\hat{\bm{z}}\right)\right]\cdot\left(\phi\,\hat{\bm{z}}\right)^*
  \mathrm{d}\Omega = 0,\\*
     \forall \bm{\phi} \in \hzcurl{\Omega}\,,\, \phi \in \hzone{\Omega}\text{.}
\end{multline}

The PML regions are incorporated naturally in the formulation by means of Equation \ref{eq:pml-mat-transf-iso},
and for simplification we have chosen to surround the domain with a Perfect Electric Conductor (PEC), resulting in
only Dirichlet boundary conditions in the electric field formulation.

Equation \ref{eq:wgma-var-final} yields in the following algebraic problem:

Find non-trivial $\left(\beta^2,{e_t}, {e_z}\right) \in (\mathbb{C} \times [\mathbb{C}]^N \times [\mathbb{C}]^M)$ such that:
\begin{equation}\label{eq:wg-discretised-form}
  \begin{bmatrix}
    A_{tt} & 0 \\
    0      & 0
  \end{bmatrix}
  \begin{Bmatrix}
    e_t \\
    e_z
  \end{Bmatrix}
  = - \beta ^2
  \begin{bmatrix}
    B_{tt} & B_{tz} \\
    B_{zt} & B_{zz}
  \end{bmatrix}
  \begin{Bmatrix}
    e_t \\
    e_z
  \end{Bmatrix}
  \text{,}
\end{equation}
where:
\addtocounter{equation}{-1}
\begin{subequations}
\begin{align}
  &\begin{multlined}\label{eq:wg-mat-1}
  [A_{tt}]_{ij} = \int_{\Omega}
  \left[
  \bm{\mu}_{r}^{-1}\left( \curlopt \bm{\varphi}_j \right) \cdot
  \left( \curlopt \bm{\varphi}_i \right)^*
\right.\\\left.
  - k_0^2 \bm{\epsilon}_r \bm{\varphi}_j\cdot \bm{\varphi}_i^*\right] d\Omega\text{,}
\end{multlined}\\
  &[B_{tt}]_{ij} = \int_{\Omega}  \left[\bm{\mu}_{r}^{-1}
    \left(\bm{\varphi}_j\times\hat{\bm{z}}\right)\right]
    \cdot
    \left(\bm{\varphi}_i\times\hat{\bm{z}}\right)^*\,
    d\Omega\text{,}\label{eq:wg-mat-2}\\
  &[B_{tz}]_{ij} = \int_{\Omega}  \left[\bm{\mu}_{r}^{-1}
    \left(\bm{\varphi}_j\times\hat{\bm{z}}\right)\right]  \cdot
  \left(\gradopt \phi_i\times\hat{\bm{z}}\right)^*\, d\Omega\text{,}\label{eq:wg-mat-3}\\
  &[B_{tz}]_{ij} = \int_{\Omega}  \left[\bm{\mu}_{r}^{-1}
    \left(\gradopt \phi_j\times\hat{\bm{z}}\right)\right]  \cdot
    \left(\bm{\varphi}_i\times\hat{\bm{z}}\right)^*\, d\Omega\text{,}\label{eq:wg-mat-4}\\
  &\begin{multlined}\label{eq:wg-mat-5}
    [B_{zz}]_{ij} = \int_{\Omega}
  \left[\bm{\mu}_r^{-1}\left(\gradopt \phi_j\times\hat{\bm{z}}\right)\right]  \cdot
  \left(\gradopt \phi_i\right)^*\\
  - k_0^2 \bm{\epsilon}_{r} \left(\phi_j\hat{\bm{z}}\right)\left( \phi_i\hat{\bm{z}}\right)^* d\Omega\text{.}
\end{multlined}
\end{align}
\end{subequations}

The generalised eigenvalue problem in Equation \ref{eq:wg-discretised-form} is Hermitian
and thus self-adjoint if the materials involved are lossless, for the tensors
$\bm{\epsilon}_r$ and $\bm{\mu}_r$ are real symmetric.
By inspection of Equation \ref{eq:pml-mat-transf} we see that it is not the case when
PMLs are involved, even if the materials were originally lossless, due to the losses
introduced by the PML region.

It is noteworthy that the De Rham compatibility between
$\hzcurl{\Omega}$ and $\hzone{\Omega}$ 
is a necessary condition\citep{demkowicz2006}
to avoid spurious modes and ensure stability.
For all the examples in this paper, the chosen spaces form an \emph{exact} sequence\citep{devloo2022},
meaning that the span of the gradient of the $\hzone{\Omega}$ space coincides
with the kernel of the curl of the $\hzcurl{\Omega}$.

\subsubsection{Scattering formulation}\label{sec:scatt}
The scattering weak form can be derived from Equation \ref{eq:helmholtz} by means of the Galerkin method.
Now, we denote by $\Omega$ the computational for the scattering analysis,
which is one dimension larger than the modal analysis domain.
The electric field $E$ is no longer restricted to the form of Equation \ref{eq:e-modal-analysis}, as most scenarios of interest involve reflections and discontinuities, and not
simply the undisturbed propagation of waveguide modes, so we simply look for $\bm{E}\in\hzcurl{\Omega}$.
The boundary is considered to be $\Gamma = \Gamma_d \cup \Gamma_m$, with $\Gamma_m$ denoting the boundary
in which the WPBC mixed boundary condition is imposed, and PEC boundaries elsewhere.

The first step is to compute the inner product with a proper test function on both sides of the equation:
\begin{multline}\label{eq:scatt-var-1}
  \int_\Omega\left\{\curlop\left[\bm{\mu}^{-1}_r\left(\curlop\bm{E}\right)\right]\right\}\cdot\bm{\phi}^*
  -k_0^2\bm{\epsilon}_ r\bm{E}\cdot\bm{\phi}^*\mathrm{d}\Omega\\
  = \int_\Omega\bm{F}\cdot\bm{\phi}^*\mathrm{d}\Omega\quad\forall \bm{\phi} \in \hzcurl{\Omega}\text{.}
\end{multline}

Using the vector identity $\divop(\bm{A}\times\bm{B})=\bm{B}\cdot(\curlop\bm{A})-\bm{A}\cdot(\curlop\bm{B})$,
\begin{multline}\label{eq:scatt-var-2}
  \int_\Omega
  \left[\bm{\mu}_{r}^{-1} \left( \curlop \bm{E}\right)\right]\cdot \left( \curlop \bm{\phi} \right)^*+\\
  \divop\left\{ \left[\bm{\mu}_{r}^{-1} \left( \curlop \bm{E}\right)\right]\times\bm{\phi}^*\right\}
  -k_0^2\bm{\epsilon}_ r\bm{E}\cdot\bm{\phi}^*\mathrm{d}\Omega\\
  = \int_\Omega\bm{F}\cdot\bm{\phi}^*\mathrm{d}\Omega\quad\forall \bm{\phi} \in \hzcurl{\Omega}\text{.}
\end{multline}

Gauss's theorem can be applied on the divergence term, element by element:

\begin{multline}\label{eq:scatt-gauss}
  \int_\Omega
  \divop\left\{ \left[\bm{\mu}_{r}^{-1} \left( \curlop \bm{E}\right)\right]\times\bm{\phi}^*\right\}
  \mathrm{d}\Omega=\\
  \sum_K \int_{\partial K}
  \left\{ \left[\bm{\mu}_{r}^{-1} \left( \curlop \bm{E}\right)\right]\times\bm{\phi}^*\right\}
  \cdot\bm{n} \mathrm{d}\Gamma=\\
  \sum_K \int_{\partial K}
  \left[\bm{\mu}_r^{-1}\curlop\bm{E}\right]\cdot\left(\bm{\phi}^*\times\bm{n}\right)\,\mathrm{d}\Gamma
\end{multline}

Assuming now a tangential continuity of $\left[\bm{\mu}_{r}^{-1} \left( \curlop \bm{E}\right)\right]$
between elements, all terms in the right-hand side of Equation \ref{eq:scatt-gauss} vanish, and
the boundary terms on $\Gamma_d$ vanish due to $\bm{\phi}\in\hzcurl{\Omega}$. Finally, the problem can be posed as:

For $\bm{F}\in\hzcurl{\Omega}$ given, find $\bm{E}\in \hcurl{\Omega}$ such that:
\begin{multline}\label{eq:scatt-var-final}
  \int_\Omega
  \left[\bm{\mu}_{r}^{-1} \left( \curlop \bm{E}\right)\right]\cdot \left( \curlop \bm{\phi} \right)^*
  - k_0^2 \bm{\epsilon}_r\bm{E}\cdot\bm{\phi}^*\mathrm{d}\Omega\\*
  +\int_{\Gamma_m}
  \left[\bm{\mu}_r^{-1}\curlop\bm{E}\right]\cdot\left(\bm{\phi}^*\times\bm{n}\right)\,\mathrm{d}\Gamma\\
  = \int_\Omega \bm{F}\cdot\bm{\phi}^*\mathrm{d}\Omega\quad\forall \bm{\phi} \in \hzcurl{\Omega}\text{.}
\end{multline}

It is worth noting that there are also no restrictions on the anisotropy of the involved media,
as long as they are still transverse anisotropic in the waveguide cross section.

\subsection{Source term modelling}

The source term $\bm{F}$ is associated with a given physical source, and in this work we discuss
two different approaches for modelling a field injection in the domain.
One alternative is to model the source term as a plane current source\citep{tsuji2002}, propagating
the incident modes away from the plane in both sides.
In Figure \ref{fig:wg-disc}, a general waveguide discontinuity is shown. The fictitious planes
$\Gamma_{\text{in}}$ and $\Gamma_{\text{out}}$ are the waveguide ports. In this approach, the waveguide ports need be PML backed, thus the extension of the domain in the $z$ direction.

\begin{figure}[ht]
  \centering
  \resizebox{\columnwidth}{!}{%
    \begin{tikzpicture}[scale=0.5]
      \draw [-{Stealth[scale=1.5]}](2.5,2.5) -- (3.75,2.5) node[right] {$z$} ;
      \draw [-{Stealth[scale=1.5]}](2.5,2.5) -- (2.5,3.75) node[above] {$x$} ;
      \draw[fill=plotcolor1,draw=none] (0,5) rectangle (6,7);
      \draw[fill=plotcolor1,draw=none] (6,4) rectangle (12,8);
      \draw[pattern=north west lines, pattern color=plotcolor2] (0,0) rectangle (12,2);
      \draw[pattern=north west lines, pattern color=plotcolor2] (0,0) rectangle (2,12);
      \draw[pattern=north west lines, pattern color=plotcolor2] (10,0) rectangle (12,12);
      \draw[pattern=north west lines, pattern color=plotcolor2] (0,10) rectangle (12,12);
      \draw [line width = 2](2,12) -- (2,0) node[pos=0.2, right] {$\Gamma_{\text{in}}$} ;
      \draw [line width = 2](10,12) -- (10,0) node[pos=0.2,left] {$\Gamma_{\text{out}}$} ;
  \end{tikzpicture}}
\caption{Computational domain representing a slab waveguide discontinuity. Dashed rectangles represent PML regions and
the propagation direction is $z$.}
  \caption{Discontinuity in a waveguide}
  \label{fig:wg-disc}
\end{figure}

To do so, the field at the waveguide source is split into the incident field and the scattered field
as $\bm{E}=\bm{E}_{\mathrm{in}}+\bm{E}_{\mathrm{scatt}}$. Assuming that only the scattered field has a continuous
tangential continuity of $\left[\bm{\mu}_{r}^{-1} \left( \curlop \bm{E}_{\textrm{scatt}}\right)\right]$,
the right-hand side term of Equation \ref{eq:scatt-gauss} vanish everywhere
except at the waveguide cross section $\Gamma$.

At $\Gamma$, we have:

\begin{equation}\label{eq:current-src-1}
  \int_\Gamma
  \llbracket
  \bm{n}\times\left[\bm{\mu}_{r}^{-1} \left( \curlop \bm{E}_{\textrm{in}}\right)\right]\rrbracket
  \cdot\bm{\phi}^* \mathrm{d}\Gamma\text{.}
\end{equation}

For computing the jump of the incident field, we write it as

\begin{subequations}
\begin{align}
  \label{eq:scatt-inc-fields}
  \begin{multlined}
  \bm{E}_{\textrm{in},1} = \sum_i\bm{E}_ie^{-j\beta_i z}\\
  =\sum_i \left(\bm{E}_{i,t}+\bm{\hat{z}}E_{i,z}\right)e^{-j\beta_i z}\text{,}
  \end{multlined}\\
  \begin{multlined}
  \bm{E}_{\textrm{in},2} = \sum_i\bm{E}'_ie^{j\beta_i z}\\
  =\sum_i \left(\bm{E}_{i,t}-\bm{\hat{z}}E_{i,z}\right)e^{j\beta_i z}\text{,}
  \end{multlined}
\end{align}
\end{subequations}
where $\bm{E}'_i$ is the backward propagating $i$-th mode.

For  waveguide modes, the curl can be computed as:
\begin{multline}
  \label{eq:curl-wg-mode-f}
  \curlop\left(\bm{E}(x,y)e^{-j\beta z}\right) =\\*
  \left[(j\beta \bm{E}_t+\gradopt E_z)\times \hat{\bm{z}}+
  \left(\frac{\partial E_y}{\partial x}-\frac{\partial Ex}{\partial y}\right)\hat{\bm{z}}\right]e^{-j\beta z}\text{,}
\end{multline}

\begin{multline}
  \label{eq:curl-wg-mode-b}
  \curlop\left(\bm{E}(x,y)e^{ j\beta z}\right) =\\*
  \left[(-j\beta \bm{E}_t-\gradopt E_z)\times \hat{\bm{z}}+
  \left(\frac{\partial E_y}{\partial x}-\frac{\partial Ex}{\partial y}\right)\hat{\bm{z}}\right]e^{j\beta z}\text{,}
\end{multline}
and therefore, the tangential part of the curl\footnote{The $z$-component of Equations \ref{eq:curl-wg-mode-f} and \ref{eq:curl-wg-mode-b} are ignored due to the cross product with $\hat{\bm{z}}$} reads as:
\begin{subequations}
\begin{align}
  \label{eq:scatt-inc-curl}
  \begin{multlined}
    \bm{n}_1\times\bm{\mu}_r^{-1}\curlop\bm{E}_{\textrm{in},1} = \\
    \bm{z}\times\left[\sum_i\left(\nabla_tE_{i,z}\times\hat{\bm{z}}\right)
    +j\beta_i \left(\bm{E}_{i,t}\times\hat{\bm{z}}\right)\right]e^{-j\beta_i z}\text{,}
  \end{multlined}\\
  \begin{multlined}
    \bm{n}_2\times\bm{\mu}_r^{-1}\curlop\bm{E}_{\textrm{in},2} = \\
    -\bm{z}\times\left[-\sum_i\left(\nabla_tE_{i,z}\times\hat{\bm{z}}\right)
    -j\beta_i \left(\bm{E}_{i,t}\times\hat{\bm{z}}\right)\right]e^{j\beta_i z}\text{,}
  \end{multlined}
\end{align}
\end{subequations}
which allows for:
\begin{multline}\label{eq:current-src-jump}
  \llbracket
  \bm{n}\times\left[\bm{\mu}_{r}^{-1} \left( \curlop \bm{E}_{\textrm{in}}\right)\right]\rrbracket=\\*
  \sum_i\left[\left(\nabla_tE_{i,z}+j\beta_i\bm{E}_{i,t}\right)\times\hat{\bm{z}}\right]2 \cos(\beta_i z)\text{,}
\end{multline}
which is the sought source term, matching the term in \citet{tsuji2002} for TE/TM modes evaluated at $z=0$.

While this approach is effective, it involves backing the port with a PML region for absorbing the
counter propagating wave, which results in an undesirable enlargement of the computational domain.
The proposed alternative is to use the WPBC, which is the subject of the following section and it allows
for the domain to be truncated exactly at the source plane.

\section{Waveguide Port Boundary Condition (WPBC)}\label{sec:wpbc}

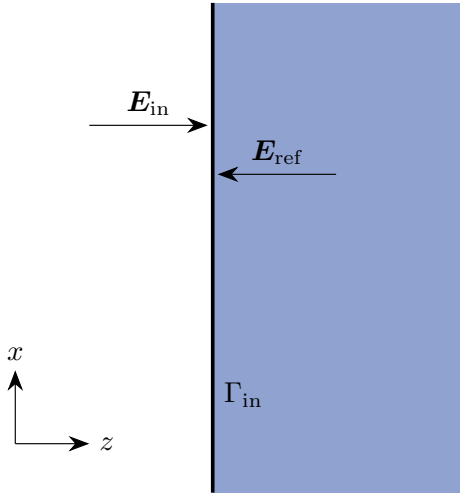
\begin{figure}[h]
  \centering
  \resizebox{0.8\columnwidth}{!}{%
    \begin{tikzpicture}[scale=6]
      \draw [-{Stealth[scale=1.5]}](0.1,0.1) -- (0.25,0.1) node[right] {$z$} ;
      \draw [-{Stealth[scale=1.5]}](0.1,0.1) -- (0.1,0.25) node[above] {$x$} ;
      \draw[fill=plotcolor1,fill opacity=0.5,draw=none] (0.5,0.0) rectangle (1.0,1.0);
      \draw[line width=1.25] (0.5,0) -- (0.5,1.0) node [pos=0.2,right]{$\Gamma_{\text{in}}$} ;
      \draw [-{Stealth[scale=1.5]}](0.75,0.65) -- (0.51,0.65) node[midway,above] {$\bm{E}_{\textrm{ref}}$} ;
      \draw [-{Stealth[scale=1.5]}](0.25,0.75) -- (0.49,0.75) node[midway,above] {$\bm{E}_{\textrm{in}}$} ;
  \end{tikzpicture}}
  \caption{Waveguide port illustration at the $-\hat{\bm{z}}$ boundary}
  \label{fig:wg-port}
\end{figure}

In Figure \ref{fig:wg-port}, an (incident) waveguide port is shown, for which we aim to
develop a boundary condition.

At $\Gamma_{\textrm{in}}$, the total field can be written as a sum of incident and reflected fields.
Assuming that the leftmost domain of $\Gamma_{\textrm{in}}$, which is not in the computational domain, consists of a waveguide,
it is reasonable to assume that the reflected fields can be decomposed as a sum of backward propagating modes
\footnote{The prime sign for backward propagating modes is omitted from now onwards, and should the
  distinction be important, it will be explicitly discussed.},
that is
\begin{equation}
  \label{eq:e-inc-ref}
  \bm{E}=\bm{E}_{\textrm{in}}+\bm{E}_{\textrm{ref}} = \bm{E}_{\textrm{in}}+\sum_i \alpha_i\bm{E}_ie^{+j\beta_i z}\text{.}
\end{equation}

This idea is closely related to the concept of eigenfunctions\citep{jin2014}:
the (possibly infinite) set of homogeneous solutions to a PDE with given boundary conditions form
a basis for the solution space of the same PDE subjected to the same boundary conditions with a given source term.

The modes of the waveguide at $\Gamma$ are a natural candidate set
for expanding the solution of the scattering problem at $\Gamma$, and that is the basis for the theory of the WPBC,
introduced in \citet{jin1994} and detailed in \citet{jin2014}. The WPBC weakly enforces that, at $\Gamma$,
the solution will be a sum of waveguide modes weighted by the appropriate coefficients $\alpha_i$ so as to
represent the reflection. In this way, the domain before $\Gamma$ is no longer discretised, as it was with
the PML backed current plane approach.

Upon inspection of Equation \ref{eq:e-inc-ref}, it is clear that with an orthogonality condition of the modes
$E_i$ at hand, an expression for $\alpha_i$ can be found.

\subsection{On the orthogonality condition of eigenmodes}\label{sec:orth}

The slab waveguide is the simplest example to illustrate the issues involving open waveguides. It consists of a plane core with higher refractive index deposited over a substrate and covered by cladding.
Due to its invariance in one of the transverse directions, the solutions for the slab waveguide can be split
into TE and TM modes. Moreover, either of the solution sets can be split into guided modes and radiation modes.

The guided modes form a discrete spectrum and are confined to the core and attenuate at infinity.
For the radiation modes, however, the situation is quite different: not only they do not attenuate at infinity,
but also they form a \emph{continuous} spectrum\citep{hu2009}.
The former implies that the solutions do not have an associated norm, much less present any sort of orthogonality.

While this behaviour might seem unusual, it is aligned with the fact that the system under analysis is a
non-conservative \emph{open} system, and thus the properties associated with Hermitian systems no longer hold\citep{sauvan2022}.

This naturally translates to the computational modelling of the slab waveguide.
As previously discussed, the computational domain must be truncated as to have finite dimensions, and this
truncation is often performed with PMLs.
As a consequence, the computational problem has materials with complex coefficients
and the resulting algebraic system is no longer Hermitian.
The continuous spectrum of radiation modes is now discretised into a discrete set of complex modes,
as thoroughly discussed in \citet{sauvan2022}, should the reader be interested in a deeper discussion on the topic.

Therefore, the orthogonality shown in \citet{jin1994} for rectangular waveguides do not hold.
While the discussion in \citet{jin2014}, allows for inhomogeneous waveguides, none of the examples present a PML or
a lossy region, and $\bm{\mu}_r$ is assumed to be constant and isotropic.

For the reasons above, we state now a general orthogonality relationship for eigenmodes, supporting losses and
anisotropy, based purely on the algebraic eigensystem of the modal analysis.

Let us denote by $(\beta_m^2,e_m)$ the solution $(\beta_m^2,\{e^m_{t},e^m_{z}\})$ from problem \ref{eq:wg-discretised-form}
associated with the $m$-th mode. From \citet[Theorem~4.30]{hanson2002}, we have then

\begin{equation}
  \label{eq:orth-1}
  (\beta_n^2-\beta_m^2)\langle Be_n, e_m^*\rangle\text{,} = (\beta_n^2-\beta_m^2)e_m^TBe_n\text{,}
\end{equation}
where $\langle\cdot,\cdot\rangle$ denotes the conjugate inner product, and we have used the fact that
the solution of the adjoint problem is the conjugate of the solution of the original problem.
Equation \ref{eq:orth-1} implies in
\begin{multline}
  \label{eq:orth-2}
  (\beta_n^2-\beta_m^2)
  \int_\Gamma
  \left\{\bm{\mu}_{r}^{-1} \left[\left( \nabla_t e_{n,z} +\bm{e}_{n,t} \right)\times\hat{\bm{z}}\right]\right\}
  \cdot\\*
  \left[\left( \nabla_t e_{m,z} +\bm{e}_{m,t} \right)\times\hat{\bm{z}}\right]\\*
  -k_0^2\left[\bm{\epsilon}_{r} \left(e_{n,z}\hat{\bm{z}}\right)\right]\cdot\left(e_{m,z}\,\hat{\bm{z}}\right)
  \mathrm{d}\Gamma = 0\text{.}
\end{multline}

The term
\begin{multline}
  \label{eq:orth-null}
  \int_\Gamma
  \left\{\bm{\mu}_{r}^{-1} \left[\left( \nabla_t e_{n,z} +\bm{e}_{n,t} \right)\times\hat{\bm{z}}\right]\right\}
  \cdot
  \left[\left( \nabla_t e_{m,z} \right)\times\hat{\bm{z}}\right]\\*
  -k_0^2\left[\bm{\epsilon}_{r} \left(e_{n,z}\hat{\bm{z}}\right)\right]\cdot\left(e_{m,z}\,\hat{\bm{z}}\right)
  \mathrm{d}\Gamma\text{,}
\end{multline}
corresponds to the right-hand side of the weak formulation evaluated at the solution $\{e^{n}_{t},e^{n}_{z}\}$
tested against the test function associated with the \emph{dofs} $\{0,e^{m}_{z}\}$.
Given that the left-hand side of the weak formulation is zero for such a test function (there are only terms
associated with $\bm{e}_{t}$), this term must also be zero and it can be eliminated, simplifying thus
Equation \ref{eq:orth-2} to 
\begin{multline}
  \label{eq:orth-3}
  (\beta_n^2-\beta_m^2)\\*
  \int_\Gamma
  \left\{\bm{\mu}_{r}^{-1} \left[\left( \nabla_t e_{n,z} +\bm{e}_{n,t} \right)\times\hat{\bm{z}}\right]\right\}
  \cdot
  \left(\bm{e}_{m,t} \times\hat{\bm{z}}\right)\\*
  \mathrm{d}\Gamma = 0\text{.}
\end{multline}

Applying the inverse transformation of Equation \ref{eq:fem-hcurl-formulation-transf}, one obtains
\begin{multline}
  \label{eq:orth-4}
  (\beta_n^2-\beta_m^2)\\*
  \int_\Gamma
  \left\{\bm{\mu}_{r}^{-1} \left[\left( j\beta_n\bm{E}_{n,t}+ \nabla_t E_{n,z} \right)\times\hat{\bm{z}}\right]\right\}
  \cdot
  \left(\bm{E}_{m,t} \times\hat{\bm{z}}\right)\\*
  \mathrm{d}\Gamma = 0\text{.}
\end{multline}

Now, we can make use of the expression for the
curl of a waveguide mode in Equation \ref{eq:curl-wg-mode-f}, repeated here for clarity
\begin{multline}
  \tag{\ref{eq:curl-wg-mode-f}}
  \curlop\left(\bm{E}(x,y)e^{-j\beta z}\right) =\\*
  \left[(j\beta \bm{E}_t+\gradop E_z)\times \hat{\bm{z}}+
  \left(\frac{\partial E_y}{\partial x}-\frac{\partial Ex}{\partial y}\right)\hat{\bm{z}}\right]e^{-j\beta z}\text{,}
\end{multline}
and, observing that the $z$-component of Equation \ref{eq:curl-wg-mode-f} does not affect the inner product,
obtain:
\begin{multline}
  \label{eq:orth-5}
  (\beta_n^2-\beta_m^2)e^{j\beta_nz}\\*
  \int_\Gamma
  \left[\bm{\mu}_{r}^{-1}\nabla\times\left(\bm{E}_ne^{-j\beta_n z}\right)\right]
  \cdot
  \left(\bm{E}_{m,t} \times\hat{\bm{z}}\right)\\*
  \mathrm{d}\Gamma = 0\text{.}
\end{multline}

Since

\begin{equation}
  \bm{H}_ne^{-j\beta_nz} = \frac{-1}{j\omega \mu_0}\bm{\mu}_r^{-1}
  \left[\nabla\times\left(\bm{E}_ne^{-j\beta_nz}\right)\right]\text{,}
\end{equation}
we have
\begin{equation}
  \label{eq:orth-semifinal}
  j\omega\mu_0(\beta_n^2-\beta_m^2)\int_\Gamma\left(\bm{H}_n\times\bm{E}_m\right)\cdot\hat{\bm{z}}\mathrm{d}\Gamma\text{.}
\end{equation}

Thus, if the modes associated with degenerate eigenvalues are orthogonalised appropriately, we can then write

\begin{equation}
  \label{eq:orth-final}
  \int_\Gamma\left(\bm{H}_n\times\bm{E}_m\right)\cdot\hat{\bm{z}}\mathrm{d}\Gamma\text{,} =
  \begin{cases}
    \kappa_m,\hfill m=n\\
    0,\hfill m\neq n
    \end{cases}\text{,}
  \end{equation}
  with
  \begin{multline}
    \label{eq:kappa-def}
    \kappa_m = \int_\Gamma\left(\bm{H}_m\times\bm{E}_m\right)\cdot\hat{\bm{z}}\mathrm{d}\Gamma =\\*
    \frac{-1}{j\omega\mu_0}
    \int_\Gamma
  \left\{\bm{\mu}_{r}^{-1} \left[\left( j\beta_m\bm{E}_{m,t}+ \nabla_t E_{m,z} \right)\times\hat{\bm{z}}\right]\right\}
  \cdot\\*
  \left(\bm{E}_{m,t} \times\hat{\bm{z}}\right)
  \mathrm{d}\Gamma
\end{multline}

\begin{figure}
\begin{subfigure}{\linewidth}
    \centering
    \includegraphics[width=.95\linewidth]{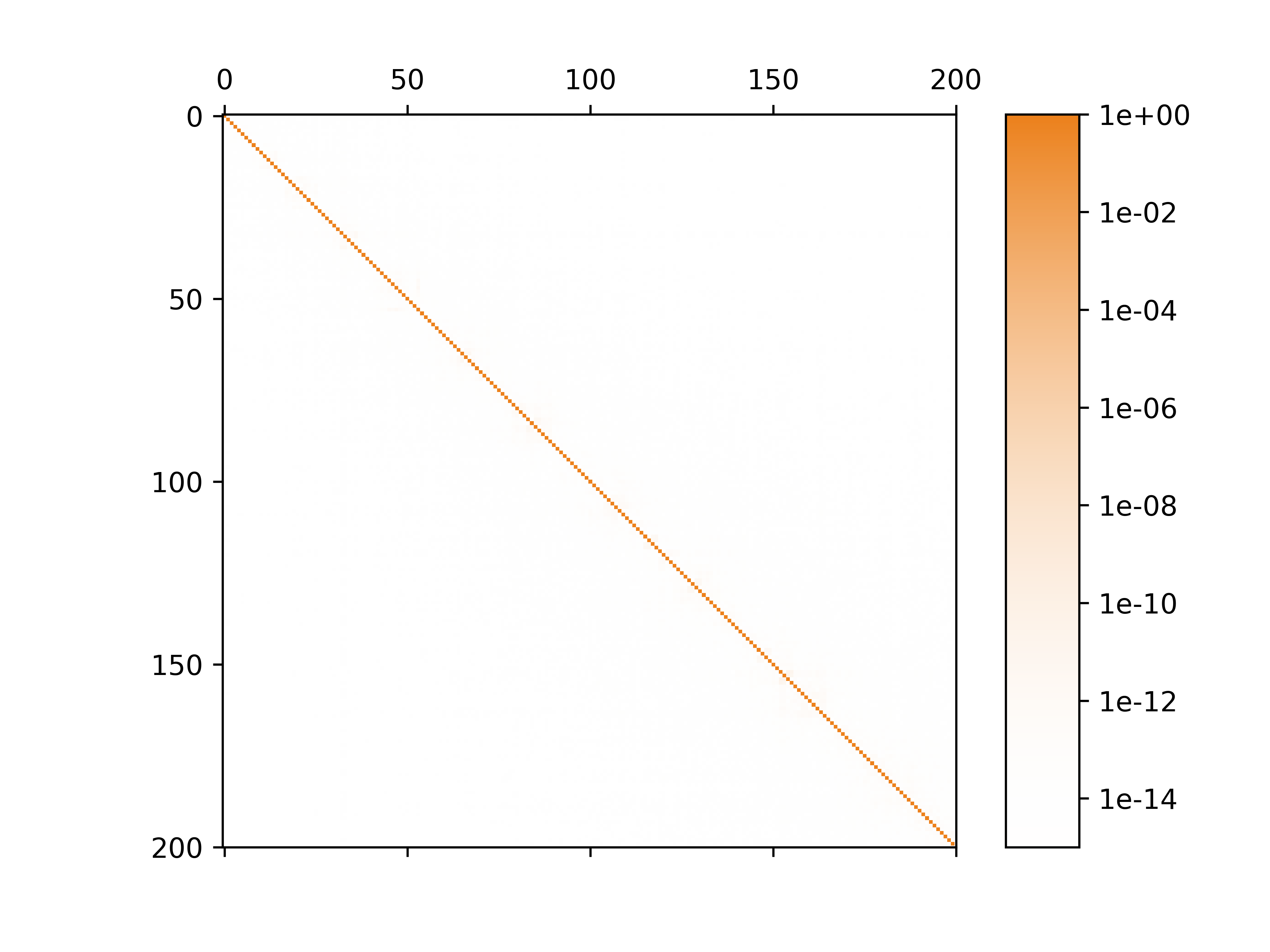}
    \caption{$\displaystyle\left\lVert\int_\Gamma \bm{H}\times\bm{E}\cdot\hat{\bm{z}}\,\mathrm{d}\Gamma\right\rVert$}
    \label{fig:orth-mat}
\end{subfigure}\\%
\begin{subfigure}{\linewidth}
    \centering
    \includegraphics[width=.95\linewidth]{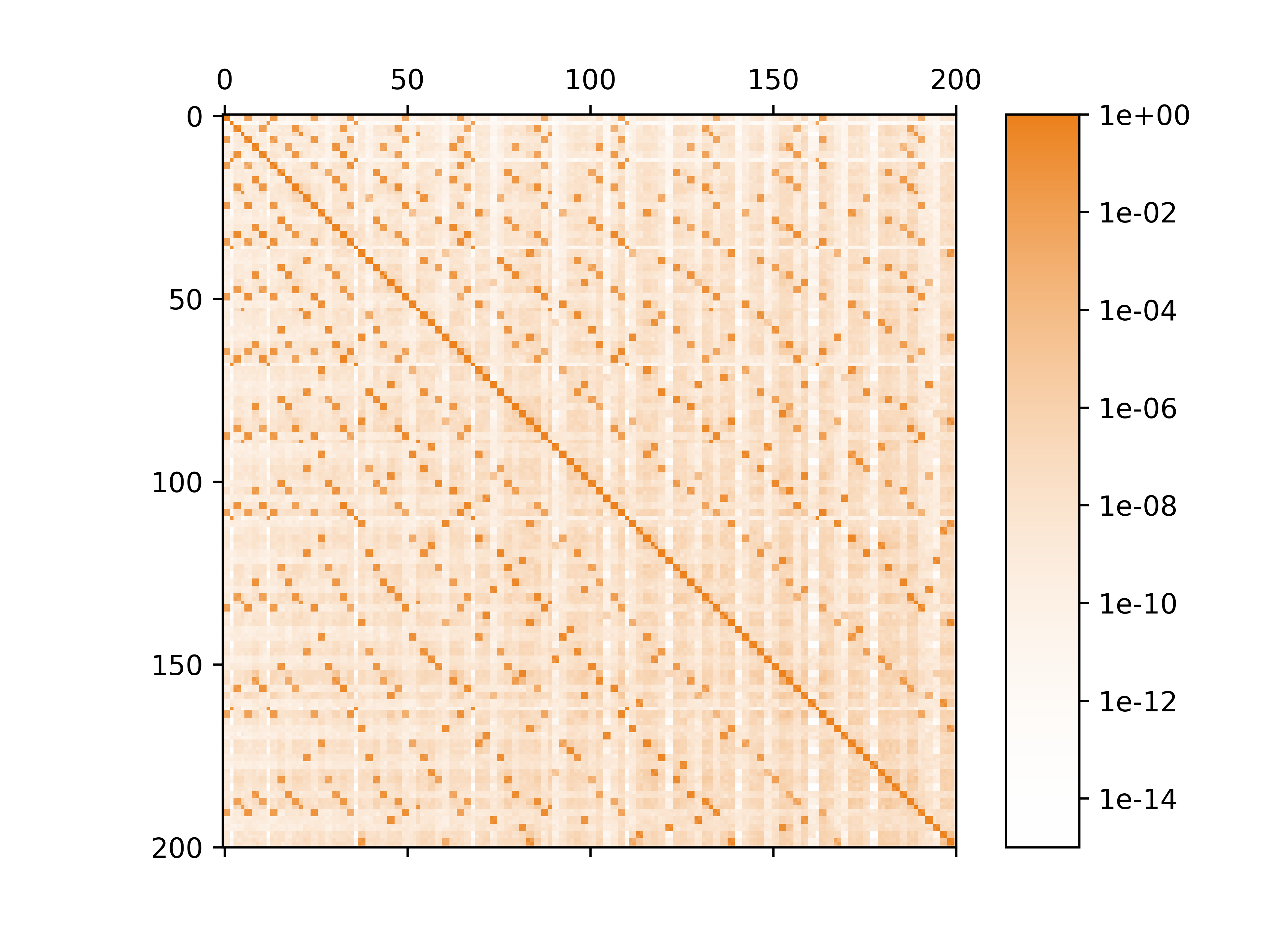}
    \caption{$\displaystyle\left\lVert\int_\Gamma \bm{H}^*\times\bm{E}\cdot\hat{\bm{z}}\,\mathrm{d}\Gamma\right\rVert$}
    \label{fig:orth-mat-conj}
\end{subfigure}
\caption{Illustrating the biorthogonality of the first 200 computed modes of a step-index optical fibre surrounded by a
  cylindrical PML. \ref{fig:orth-mat} shows the magnitude of the integrated non-conjugate cross product of the electric and magnetic fields.
  In \ref{fig:orth-mat-conj} one can see that the conjugate biorthogonality does not hold if the system is non-Hermitian.}
\label{fig:orth-mat-both}
\end{figure}

The results of this section are summarised in Figure \ref{fig:orth-mat-both},
where the magnitude of the tangential component of the cross product between the magnetic and electric fields
is plotted for the first two hundred computed modes of a step-index optical fibre, illustrating that for non-Hermitian
systems, the orthogonality holds only for the non-conjugated cross product.
\subsection{Waveguide port boundary condition}

With the orthogonality condition at hand, we can finally derive the WPBC.
The derivation will be performed assuming that the port is aligned with the $xy$ plane
and the outward normal vector is $\hat{\bm{z}}$ in the negative direction, as illustrated
in Figure \ref{fig:wg-port}.

Repeating for clarity Equation \ref{eq:e-inc-ref}:
\begin{equation}
  \tag{\ref{eq:e-inc-ref}}
  \bm{E}=\bm{E}_{\textrm{in}}+\bm{E}_{\textrm{ref}} = \bm{E}_{\textrm{in}}+\sum_i \alpha_i\bm{E}_ie^{+j\beta_i z}\text{,}
\end{equation}
it should be clear that Equations \ref{eq:orth-final} and \ref{eq:kappa-def} allows for a closed expression for $\alpha_i$, which is obtained by taking the cross product of \ref{eq:e-inc-ref} with $\bm{H}_m$ followed by the dot product with
the outward normal vector $-\bm{z}$ and integrating over the boundary $\Gamma$, resulting in

\begin{equation}
  \label{eq:wpbc-alpha-1}
  \alpha_m = \frac{e^{-j\beta_m z}}{\kappa_m}
  \int_\Gamma \bm{H}_m\times\left(\bm{E}-\bm{E}_{\textrm{in}}\right)\cdot\hat{\bm{z}}\,\mathrm{d}\Gamma\text{.}
\end{equation}

Since our formulation involves only the $\bm{E}$ field, we can eliminate $\bm{H}_m$ through Faraday's Law
applied to a backward propagating waveguide mode, that is

\begin{equation}
  \label{eq:faraday-wg-minus-z}
  \bm{H}_m e^{j\beta_m z} = \frac{-1}{j\omega \mu_0}\bm{\mu}_r^{-1}\curlop\left(\bm{E}_m e^{j\beta_m z}\right)\text{,}
\end{equation}
which, if applied in \ref{eq:wpbc-alpha-1}, results in 
\begin{multline}
  \label{eq:wpbc-alpha-semifinal}
  \alpha_m = 
  \frac{-e^{-j\beta_m z}}{j\omega\mu_0\kappa_m}\\
  \int_\Gamma \left[\bm{\mu}_r^{-1}\curlop\left(\bm{E}_m e^{j\beta_m z}\right)e^{-j\beta_m z}\right]\times\left(\bm{E}-\bm{E}_{\textrm{in}}\right)\cdot\hat{\bm{z}}\,\mathrm{d}\Gamma\text{.}
\end{multline}

Upon the expansion of the expression for the curl of a backward propagating waveguide mode,
\begin{multline}
  \label{eq:wpbc-alpha-final}
  \alpha_m = \frac{-e^{-j\beta_m z}}{j\omega\mu_0\kappa_m}\\*
  \int_\Gamma \bm{\mu}_r^{-1}\left[\left(-j\beta_m\bm{E}_{m,t}+\gradopt E_{m,z}\right)\times\bm{z}\right]\\*
  \cdot\left[\left(\bm{E}-\bm{E}_{\textrm{in}}\right)\times\hat{\bm{z}}\right]\,\mathrm{d}\Gamma\text{.}
\end{multline}

Equation \ref{eq:wpbc-alpha-final} is a closed expression for $\alpha_m$ in terms of the
source field $\bm{E}_{\textrm{in}}$ and the total field $\bm{E}$,
therefore suitable for usage in a mixed boundary condition term
for the scattering formulation. Since the boundary term of Equation \ref{eq:scatt-var-final} is of the form
\begin{equation}
  \label{eq:scatt-bc-term}
  \int_\Gamma\left[\bm{\mu}_r^{-1}\curlop\bm{E}\right]\cdot\left(\bm{\phi}^*\times\bm{n}\right)\,\mathrm{d}\Gamma\text{,}
\end{equation}
we manipulate Equation \ref{eq:e-inc-ref} as to obtain an expression for $\bm{\mu}_r^{-1}\curlop\bm{E}$.
Taking its curl and applying the inverse relative permeability tensor:
\begin{multline}
  \bm{\mu}_r^{-1}\curlop\bm{E} = \\
  \bm{\mu}_r^{-1}\curlop\bm{E}_{\textrm{in}}+\sum_i\alpha_i \bm{\mu}_r^{-1}
  \left[\curlop \left(\bm{E}_i e^{j\beta_i z}\right)\right]=\\
  \bm{\mu}_r^{-1}\curlop\bm{E}_{\textrm{in}}+\\
  \sum_i\alpha_i \bm{\mu}_r^{-1}
  \left[\left(-j\beta_i\bm{E}_{i,t}+\gradopt E_{i,z}\right)\times\bm{z}\right]e^{j\beta_i z}\text{,}
\end{multline}
and by application of the expression for $\alpha_i$ given in Equation \ref{eq:wpbc-alpha-final}:
\begin{multline}
  \bm{\mu}_r^{-1}\curlop\bm{E} = \bm{\mu}_r^{-1}\curlop\bm{E}_{\textrm{in}}\\*
  -\sum_i\frac{1}{j\omega\mu_0\kappa_i}\left[\left(-j\beta_i \bm{E}_{t,i}+\gradopt E_{z,i}\right)\times\bm{z}\right]\\*
  \int_\Gamma \bm{\mu}_r^{-1}\left[\left(-j\beta_i\bm{E}_{,t}+\gradopt E_{i,z}\right)\times\bm{z}\right]\cdot\left[\left(\bm{E}-\bm{E}_{\textrm{in}}\right)\times\hat{\bm{z}}\right]\,\mathrm{d}\Gamma\text{.}
\end{multline}

The curl at the boundary can thus be expressed as
\begin{equation}
  \label{eq:scatt-bc-mixed-1}
  \bm{\mu}_r^{-1}\curlop\bm{E}+\bm{P}(\bm{E}) = \bm{U}_{\textrm{in}}\text{,}
\end{equation}
where
\begin{multline}
  \label{eq:p-e-plus}
  P(\bm{E}) = \\*
  \sum_{i}\frac{1}{j\omega\mu_0\kappa_i}\bm{\mu}_r^{-1}
      \left[\left(-j\beta_i \bm{E}_{t,i}+\gradopt E_{z,i}\right)\times\bm{z}\right]\\*
      \int_{\Gamma}\bm{\mu}_r^{-1}
      \left[\left(-j\beta_i \bm{E}_{t,i}+\gradopt  E_{z,i}\right)\times\bm{z}\right]\cdot\left(\bm{E}\times\bm{z}\right)\,\mathrm{d}\Gamma
\end{multline}
and
\begin{equation}\label{eq:u-inc}
  \bm{U}^{\textrm{in}} = \bm{\mu}_r^{-1}\curlop\bm{E}^{\textrm{in}}+P(\bm{E}^{\textrm{in}})\text{.}
\end{equation}
resulting in
\begin{equation}
  \label{eq:scatt-bc-mixed-final}
  \int_\Gamma\bm{P}(\bm{E})\cdot\left(\bm{\phi}^*\times\bm{z}\right)\,\mathrm{d}\Gamma =
    \int_\Gamma\bm{U}_{\textrm{in}}\cdot\left(\bm{\phi}^*\times\bm{z}\right)\,\mathrm{d}\Gamma\text{,}
\end{equation}
since at this boundary $\bm{n}=-\bm{z}$. Equations \ref{eq:scatt-bc-mixed-1} to \ref{eq:scatt-bc-mixed-final}
are similar to what has been developed in the literature
for the isotropic (and with homogeneous permeability) case.
Upon inspection of Equation \ref{eq:scatt-bc-mixed-final}, one realises that the double integrals mean that $P(\bm{E})$ will
need to be evaluated for each trial function at $\Gamma$, which most likely translates to cumbersome implementation issues.

At this point, it might be useful to step back and reflect on the main aspects of the WPBC.
First, it enforces that the solution at $\Gamma$ is a linear combination of waveguide modes.
Moreover, it computes the coefficients $\alpha_m$ such that the resulting combination represents the reflection
at the boundary. Both aspects are enforced \emph{weakly}, since they are imposed by the weak formulation and
not at the approximation space level. It turns out that the WPBC can be greatly simplified if one of these
aspects is strongly enforced through restriction of the approximation space, which is the subject of the next
section.
\subsection{Approximation space restriction}\label{sec:approx-space-restr}

The inspection of Equations \ref{eq:p-e-plus} and \ref{eq:orth-final} provides
an insight on how the WPBC could be simplified should the trial functions at the boundary
coincide with the tangential field components $\bm{E}_{j,t}$.
Assuming that all the modes at the boundary were computed, this operation is equivalent
to a change of basis or, when just a few modes are applied, to apply a reduced basis
for the degrees of freedom at the port.
 
A comprehensive description of the implementation and the intricacies of the procedure is the subject of ongoing work; however, an illustrative example is provided here.

\begin{figure}[ht]
  \centering
  \begin{subfigure}{0.5\linewidth}
    \begin{tikzpicture}[scale=3.5]
      \draw[draw,fill=plotcolor1,opacity=0.8] (0,0)--(1,0)--(0,1);
      \draw[draw,fill=plotcolor2,opacity=0.8] (1,0)--(0.5,0.5)--(1,1);
      \draw[draw,fill=plotcolor2,opacity=0.8] (0,1)--(0.5,0.5)--(1,1);
      \draw (0.3, 0.3) node[above] {$K_1$};
      \draw (0.5, 0.75) node[above] {$K_2$};
      \draw (0.75, 0.5) node[right] {$K_3$};
    \end{tikzpicture}
  \end{subfigure}\begin{subfigure}{0.45\linewidth}
      \begin{tikzpicture}[scale=1.6]

      \node at (-1,0)[circle,fill,inner sep=1pt]{};
      \node at ( 1,0)[circle,fill,inner sep=1pt]{};

      \node at (-1,-2)[circle,fill,inner sep=1pt]{};
      \node at ( 0,-2)[circle,fill,inner sep=1pt]{};
      \node at ( 1,-2)[circle,fill,inner sep=1pt]{};

      \draw[-] (-1, 0) -- (1, 0) node[right] {};

      \draw (-1, 1) node[above] {$\phi_1$};
      \draw (1, 1) node[above] {$\phi_2$};
      \draw[domain=-1:1, smooth, variable=\x, plotcolor1] plot ({\x}, {(1-\x)/2});
      \draw[domain=-1:1, smooth, variable=\x, plotcolor1] plot ({\x}, {(1+\x)/2});

      \draw[-] (-1, -2) -- (1, -2) node[right] {};

      \draw (-1,-1) node[above] {$\phi_3$};
      \draw (0,-1) node[above] {$\phi_4$};
      \draw (1,-1) node[above] {$\phi_5$};
      \draw[domain= 0:1, smooth, variable=\x, plotcolor2] plot ({\x}, {1-\x-2});
      \draw[domain= 0:1, smooth, variable=\x, plotcolor2] plot ({\x}, {\x-2});
      \draw[domain= -1:0, smooth, variable=\x, plotcolor2] plot ({\x}, {-\x-2});
      \draw[domain= -1:0, smooth, variable=\x, plotcolor2] plot ({\x}, {\x-1});
    \end{tikzpicture}
  \end{subfigure}
  \caption{Left: A non-conforming $\hone{\Omega}$ mesh presenting a hanging node. Right: The traces
  of the basis functions that have non-vanishing traces on hypotenuse of the blue triangular element.}
  \label{fig:hanging-node}
\end{figure}
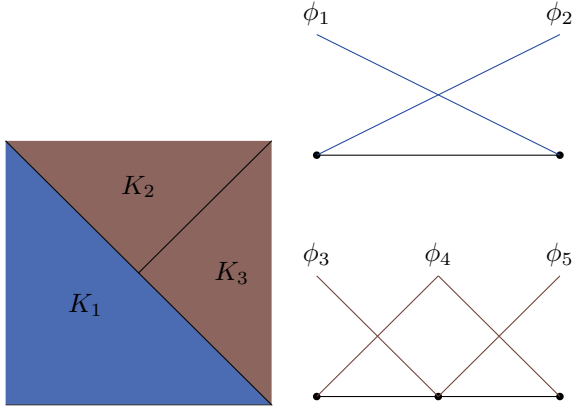

Approximation space restriction has been present in some - but not all - FEM codes for a couple of decades now,
specially in the context of \emph{hanging nodes} support, as shown in Figure \ref{fig:hanging-node}.

In this example, an $\hone{\Omega}$ space is used for simplicity, and one can see
that the mesh is not $\hone{\Omega}$-conforming as shown by the relevant traces.
One approach to make this mesh conforming without further \emph{h}-refinement is to \emph{restrict}
the functions $\phi_3$, $\phi_4$ and $\phi_5$, so they are no long free \emph{dofs}.

In this case, we can see that defining basis functions $\phi_3'$ and $\phi_4'$ such that
\begin{align}
  \label{eq:hanging-node-eqs}
  \phi_3' = \phi_3+\frac{1}{2}\phi_4\\
  \phi_4' = \phi_5+\frac{1}{2}\phi_4
\end{align}
will result in a conforming mesh. Therefore, we say that the functions $\{\phi_3,\phi_4,\phi_4\}$ are \emph{restricted}
to $\phi_3'$ and $\phi_4'$, as Equation \ref{eq:hanging-node-eqs} must always
hold for the space to be $\hone{\Omega}$-conforming.

In terms of implementation, the dependency coefficients can be computed through an $\ltwo{\partial K}$ projection
of the traces and, naturally, this dependency must be stored in the computational mesh.
Then, the element matrix is computed as if no dependencies were present,
and subsequently the dependencies are computed before contributing to the global matrix.

Assuming that the functions $\phi_3$, $\phi_4$ and $\phi_6$ have local support in the
element $K_2$, this element will contribute to the equations of $\phi_3'$, $\phi_4'$ and $\phi_6$ and its dependency matrix reads as

\begin{equation}
  \bm{D}_2 = \begin{pmatrix}
    1 & 0 & 0 \\
    0 & \frac{1}{2} & 0 \\
    0 & 0  & 1
  \end{pmatrix}\text{,}
\end{equation}
and the element matrix is transformed as
\begin{equation}
  \bm{K}_2' = \bm{D}_2^\dagger \bm{K}_2 \bm{D}_2
\end{equation}
as to incorporate the restrictions, where $^\dagger$ stands for conjugate transpose.

Similarly, for the $K_3$ element,
assuming that its local-supported functions are ordered as $\phi_4$, $\phi_5$ and $\phi_7$ results in the following dependency matrix

\begin{equation}
  \bm{D}_3 = \begin{pmatrix}
    \frac{1}{2} & 0 & 0 \\
    0 & 1 & 0 \\
    0 & 0  & 1
  \end{pmatrix}\text{.}
\end{equation}

\subsection{Waveguide port boundary condition with approximation space restriction}
The procedure described in Subsection \ref{sec:approx-space-restr} can be applied to the WPBC: if
one restricts, at $\Gamma$, the approximation space to the tangential components of the waveguide modes
$\{\bm{E}_{t,i}\}$, then the weak formulation can be simplified, as it is now only responsible
for computing the coefficients $\alpha_i$. The simplification becomes clear when we evaluate $P(\bm{E})$
for a given transverse modal component, that is:

\begin{multline}
  \label{eq:p-e-plus-restricted-semifinal}
  P(\bm{E}_{t,j}) = \sum_{i}\bm{\mu}_r^{-1}
      \left[\left(-j\beta_i \bm{E}_{t,i}+\gradopt E_{z,i}\right)\times\bm{z}\right]\frac{1}{\kappa_i}\\
      \underbrace{
        \begin{multlined}\frac{1}{j\omega\mu_0}\\
          \int_{\Gamma}\bm{\mu}_r^{-1}
        \left[\left(-j\beta_i \bm{E}_{t,i}+\gradopt  E_{z,i}\right)\times\bm{z}\right]\cdot\left(\bm{E}_{t,j}\times\bm{z}\right)\,\mathrm{d}\Gamma\text{,}\end{multlined}}_{-\kappa_i \delta_{ij}}
  \end{multline}
resulting in
  \begin{multline}
  \label{eq:p-e-plus-restricted-final}
  P(\bm{E}_{t,j}) = -\bm{\mu}_r^{-1}
      \left[\left(-j\beta_j \bm{E}_{t,j}+\gradopt E_{z,j}\right)\times\bm{z}\right]\text{.}
\end{multline}

The biorthogonality was responsible for eliminating the double integral of the boundary term. We can apply the
same principle to the source term $\bm{U}^{\textrm{in}}$: let us assume now that the incident wave can be projected
over forward propagating waveguide modes, $\bm{E}^{\textrm{in}} = \sum_k\alpha_k \bm{E}_k'e^{-j\beta_k z}$, where
$\bm{E}_k' = \bm{E}_{k,t}-\hat{\bm{z}}E_{z,t}$. This results in

\begin{multline}
  \label{eq:u-inc-restr}
  \bm{U}^{\textrm{in}} = \sum_k\left\{
    \bm{\mu}_r^{-1}\left[\left(j\beta_k\bm{E}_{k,t}-\gradopt E_{k,z}\right)\times\hat{\bm{z}}\right]-\right.\\
  \left.
    \bm{\mu}_r^{-1}\left[\left(-j\beta_k\bm{E}_{k,t}+\gradopt E_{k,z}\right)\times\hat{\bm{z}}\right]
  \right\}\alpha_ke^{-j\beta_k z}\text{,}
\end{multline}
and
\begin{multline}
  \label{eq:u-inc-restr-2}
  \bm{U}^{\textrm{in}} = -2\sum_k
  \bm{\mu}_r^{-1}\left[\left(-j\beta_k\bm{E}_{k,t}+\gradopt E_{k,z}\right)\times\hat{\bm{z}}\right]\\
  \alpha_ke^{-j\beta_k z}\text{.}
\end{multline}

Finally, we can now write the final expressions for the boundary term at the weak formulation:

\medskip
\noindent\underline{$+\hat{\bm{z}}$ boundary:}

\noindent
\begin{equation}\label{eq:wpbc-plus-z}
\boxed{
  \begin{aligned}
    \int_\Gamma
    \left\{\bm{\mu}_r^{-1}\left[\left(j\beta_j\bm{E}_{j,t}+\gradopt E_{j,z}\right)\times\hat{\bm{z}}\right]\right\}
    \cdot \left(\bm{E}_{i,t}^*\times\hat{\bm{z}}\right)\,\mathrm{d}\Gamma=\\
    2\sum_k\alpha_k e^{ j\beta_k z}\\
    \int_\Gamma\left\{\bm{\mu}_r^{-1}\left[\left(j\beta_k\bm{E}_{k,t}+\gradopt E_{k,z}\right)\times\hat{\bm{z}}\right]\right\}\cdot \left(\bm{E}_{i,t}^*\times\hat{\bm{z}}\right)\,\mathrm{d}\Gamma
  \end{aligned}}
\end{equation}

\noindent\underline{$-\hat{\bm{z}}$ boundary:}

\noindent
\begin{equation}\label{eq:wpbc-minus-z}
\boxed{
  \begin{aligned}
   \int_\Gamma
    \left\{\bm{\mu}_r^{-1}\left[\left(j\beta_j\bm{E}_{j,t}-\gradopt E_{j,z}\right)\times\hat{\bm{z}}\right]\right\}
    \cdot \left(\bm{E}_{i,t}^*\times\hat{\bm{z}}\right)\,\mathrm{d}\Gamma=\\
    2\sum_k\alpha_k e^{-j\beta_k z}\\
    \int_\Gamma\left\{\bm{\mu}_r^{-1}\left[\left(j\beta_k\bm{E}_{k,t}-\gradopt E_{k,z}\right)\times\hat{\bm{z}}\right]\right\}\cdot \left(\bm{E}_{i,t}^*\times\hat{\bm{z}}\right)\,\mathrm{d}\Gamma
  \end{aligned}}
\end{equation}

Inspecting the final expressions should suffice to illustrate how the WPBC is simplified when restricting
the approximation space. As a final remark, it is important to assure that the waveguide modes used in
the expressions correspond to the boundary orientation. However, this is quite simple to deal with,
since changing the sign of the $\hat{\bm{z}}$-component suffices
for alternating between forward/backward propagating modes.

\subsection{Implementation aspects}

The contributions of the WPBC can be computed in a post-processing manner, using the solutions of the modal analysis.
Thus, the FEM scheme can simply assemble the matrix associated with Equation \ref{eq:scatt-var-final}
and subsequently add the WPBC contributions.
Inspecting Equations \ref{eq:wpbc-plus-z} and \ref{eq:wpbc-minus-z} one can observe that, in the Hermitian problems,
the BC terms result in a diagonal matrix. For the more general case, this doesn't hold: the WPBC contributions
is actually a structurally symmetric complex matrix. Therefore, the resulting algebraic system is now no longer
complex symmetric.

One interesting remark is that, due to the approximation space restriction, at the waveguide port, the scheme will have
one \emph{dof} per waveguide mode. Therefore, if the modes are normalised according to a quantity of interest,
such as power, then the \emph{dof} itself provides automatically a result that normally one would need to
post-process to obtain, \emph{e.g.}, the power carried by each mode.

Before discussing the results, one important remark is needed on an efficient implementation of the WPBC: each element adjacent to the WPBC will contribute to all the equations
associated with the modes in this boundary. Therefore, the element matrix obtained after
the approximation space restriction can become quite large when using hundreds of modes
in the corresponding port.
Should the computational domain present reduced dimensions in the propagation direction,
which is often desirable in the analysis of a waveguide discontinuity, one can expect
to an impact in the assembly time due to the elements adjacent to the waveguide ports
when employing a large number of modes at the WPBC.

In order to mitigate this impact, these elements were grouped in patches, with one
element matrix by patch. The effect on the sparsity pattern of the resulting matrix
was not found to be significant, but it should be straightforward to recover the matrix
with the original sparsity pattern, without the element groups.

\section{Numerical results}\label{sec:res}

In this section, numerical results are presented as to proceed with the discussion on the WPBC.
All the results were obtained using the \texttt{NeoPZ}
library \citep{neopzgit}, by means of the \texttt{wgma} library\citep{wgmagit}, which acts as a \emph{wrapper}
for electromagnetic problems. All the presented examples are present in the \texttt{wgma} repository,
and the corresponding scripts will be referenced accordingly.
All meshes were generated using the \texttt{gmsh}\citep{geuzaine2009} python interface.
The Intel MKL library was used both for the sparse direct solver as for sparse BLAS-2 level operations,
and the \texttt{SLEPc} library\citep{slepc05} was used for its Krylov-based eigensolver.
Finally, graphic output was exported in \texttt{vtk} format and visualised in \texttt{ParaView}\citep{paraview2015}.

\subsection{Initial validation: effect of unstructured meshes}
In order to validate the scheme, a simple slab waveguide is analysed.
Given the planar nature of the waveguide, the modal and scattering analyses from \citet{tsuji2002} are performed,
using $\hone{\Omega}$-conforming elements in one and two dimensions, respectively.

The first experiment consists of a short section of a symmetric slab with width $w=\qty{1}{\micro\metre}$, truncated by WPBC
boundaries denoted by $\Gamma_{\mathrm{in}}$ and $\Gamma_{\mathrm{out}}$.
At $\Gamma_{\mathrm{in}}$, a combination of the guided modes of the slab, obtained through the modal analysis,
is injected into the domain. As there is no discontinuity, the incident field should undisturbingly propagate along the
waveguide, presenting no reflections.
Denoting the incident field by $\bm{E}_{\mathrm{inc}} = \sum_k \alpha_k \bm{E}_k'e^{-j\beta_k z}$, at a distance $d$ from
$\Gamma_{\mathrm{in}}$, the solution is expected to be

\begin{equation}
  \label{eq:slab-val-1}
  \bm{E}_{\mathrm{ref}} = \sum_k \alpha_k \bm{E}_k'e^{-j\beta_k d}\text{.}
\end{equation}

\begin{figure}[ht]
  \centering
  \begin{tikzpicture}
    \node[anchor=south west, inner sep=0] (X) at (0,0){\includegraphics[width=1\linewidth]{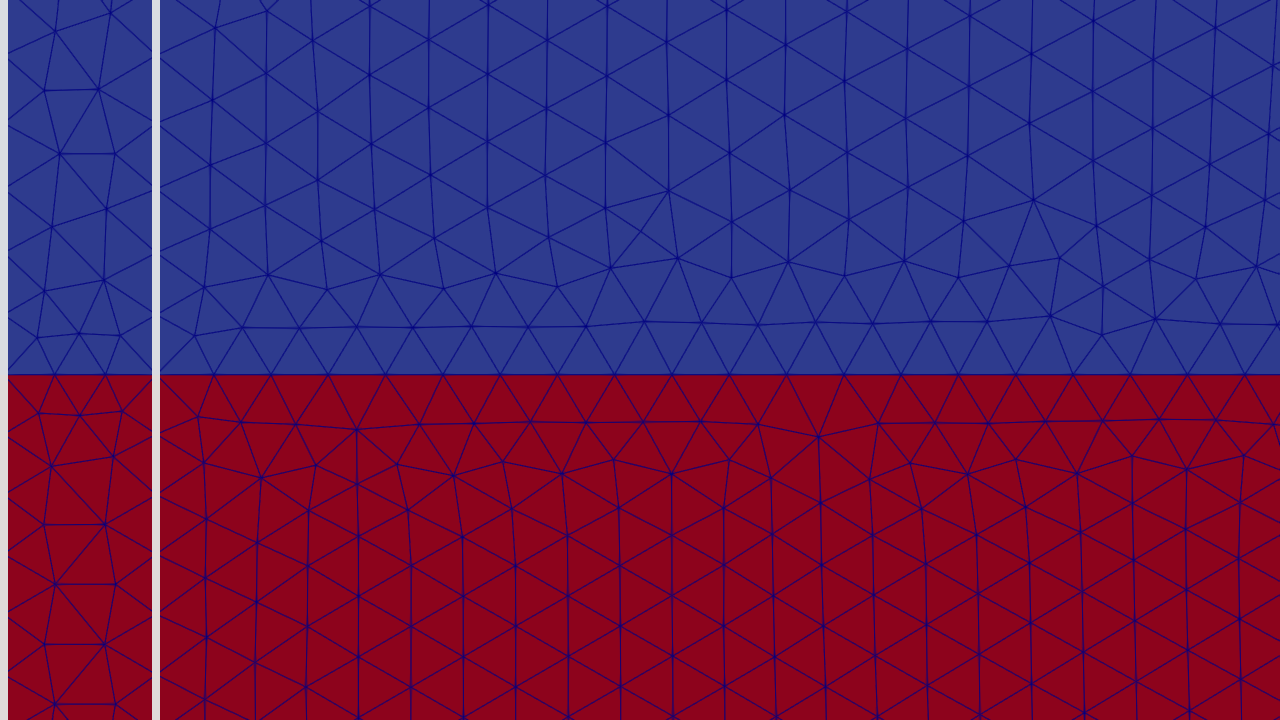}};%
  \begin{scope}[x={(X.south east)},y={(X.north west)}]%
    \draw[fill = white,fill opacity=0.8, text opacity=1] (0.01,0.05) rectangle (0.08,0.15)  node[pos=.5]{$\Gamma_{\mathrm{in}}$};
		\draw[fill = white,fill opacity=0.8, text opacity=1] (0.13,0.05) rectangle (0.20,0.15)  node[pos=.5]{$\Gamma_{\mathrm{in}}^{e}$};
		\end{scope}%
  \end{tikzpicture}
  \caption{Detail on the slab waveguide mesh used in the validation of the WPBC. The upper part represents the cladding,
  the lower part the core. Thick lines represents the unidimensional subdomains $\Gamma_{\mathrm{in}}$, the waveguide port, and $\Gamma_{\mathrm{in}}^{e}$, where the solution is evaluated.}
  \label{fig:slab-validation}
\end{figure}

As to verify the expected solution, the meshing of the slab was performed such that there is a plane (line) $\Gamma_{\mathrm{in}}^e$ at a distance $d$ from $\Gamma_{\mathrm{in}}$ with nodes matching the nodes in $\Gamma_{\mathrm{in}}$,
as shown in Figure \ref{fig:slab-validation}.

Then, one can compare the difference between $\bm{E}_{\mathrm{WPBC}}$,
the field obtained through the scattering analysis, and $\bm{E}_{\mathrm{ref}}$. Moreover, $\bm{E}_{\mathrm{ref}}$
can also be compared with $\bm{E}_{\mathrm{pml}}$,
the field obtained in the scattering analysis with PML-backed plane sources.

The simulation was performed at the operational wavelength $\lambda=\qty{1.55}{\micro\metre}$
and the PML parameters were computed according to Subsection \ref{sec:pml-param}
resulting in $s_x=39.76$ and $s_z = 59.64$. The used polynomial order for the $\hone{\Omega}$ elements was 4,
and the remaining parameters are available at the script \texttt{py/scripts/slab\_disc\_validation.py}.

\begin{figure}[ht]
  \centering
  \includegraphics[width=\linewidth]{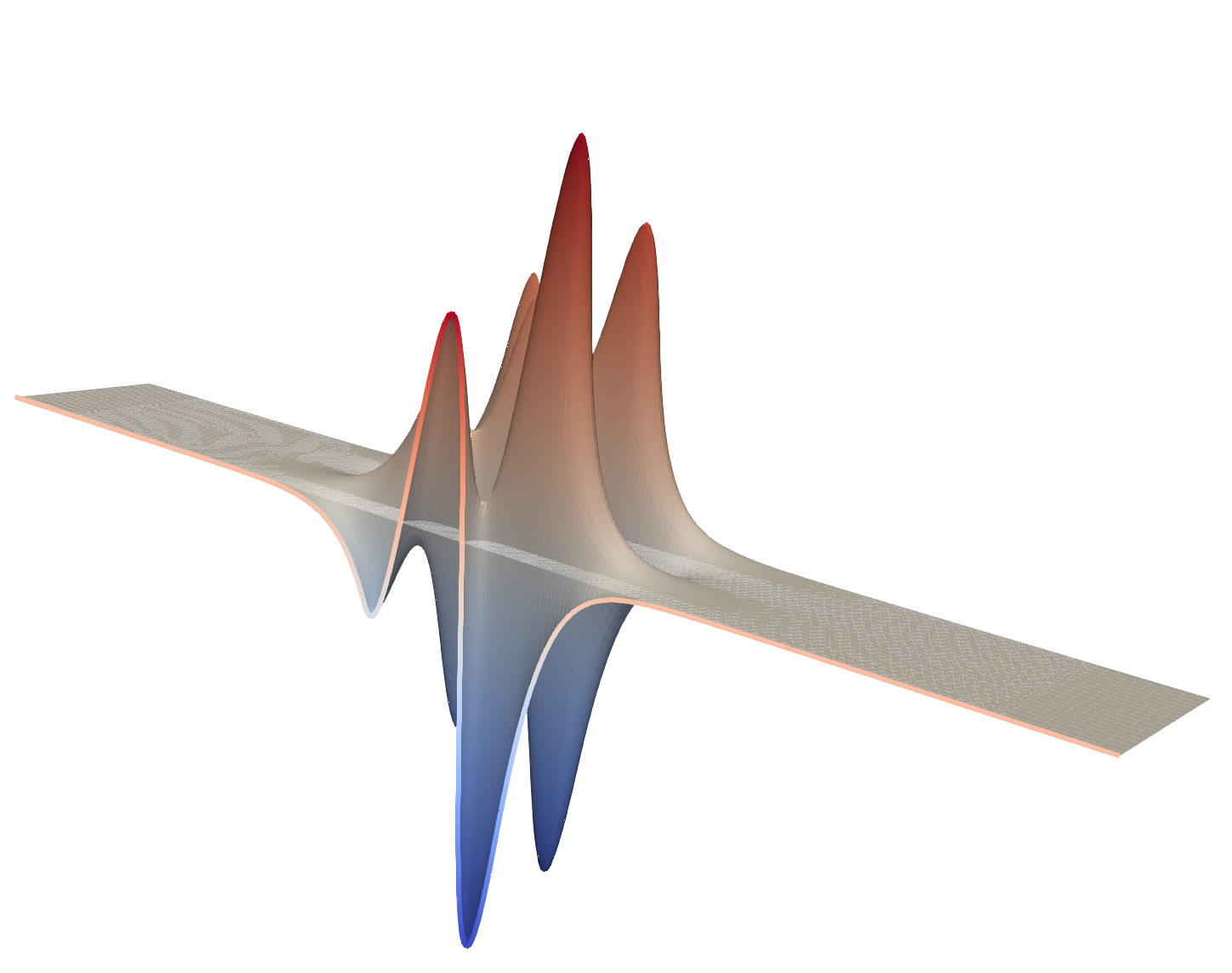}
  \caption{Combination of guided modes propagating in a slab waveguide. Results obtained using the WPBC boundary condition on both domain truncations. The incident field is represented through a thicker line.}
  \label{fig:slab-val-plot}
\end{figure}

For $\bm{E}_{\mathrm{inc}}$ set as $\alpha_1 = 0.5$, $\alpha_2=2$ and $\alpha_3 = 2.5$, the computed field is shown
in Figure \ref{fig:slab-val-plot} and the obtained relative errors are

\begin{align}
\begin{split}
  \label{eq:slab-val-res}
  \frac{\lVert\bm{E}_{\mathrm{ref}}-\bm{E}_{\mathrm{PML}}\rVert}{\lVert\bm{E}_{\mathrm{ref}}\rVert} = \num{2.235e-06} \\
  \frac{\lVert\bm{E}_{\mathrm{ref}}-\bm{E}_{\mathrm{WPBC}}\rVert}{\lVert\bm{E}_{\mathrm{ref}}\rVert} = \num{6.203e-08} \text{.}
\end{split}
\end{align}

The WPBC approach resulted in an error almost $10^2$ smaller than the one in the PML-backed propagation,
which is known to introduce numerical reflections\citep{johnson2008}. Moreover, Equation \ref{eq:slab-val-1} provides an insight on how to determine the number of modes
to be used in the WPBC, which is the subject of the next section.

\begin{figure}[ht]
  \centering
  \includegraphics[width=\linewidth]{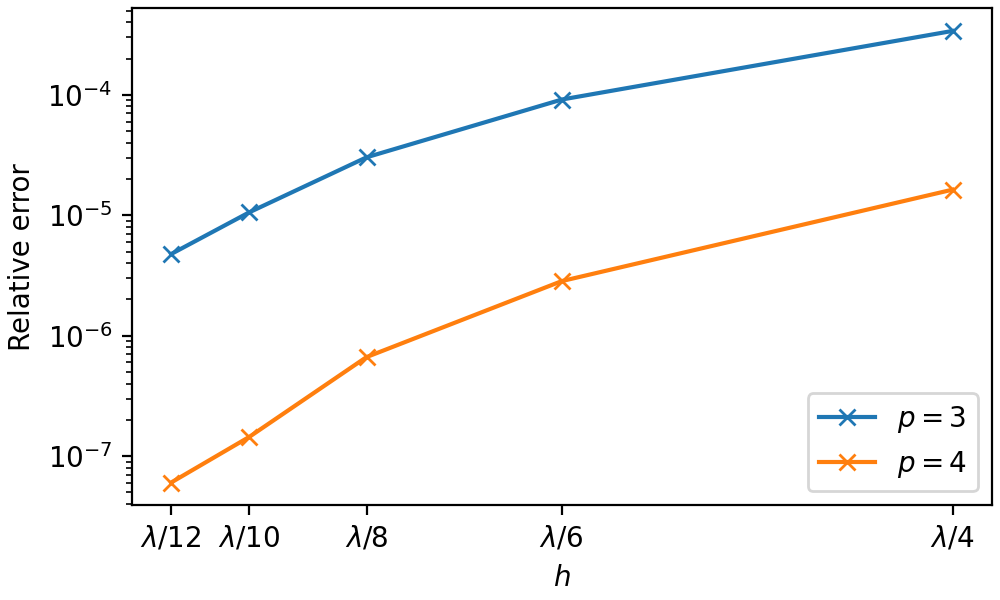}
  \caption{Relative error between WPBC approximation and analytical mode propagation
    as a function of element size and polynomial order for the slab waveguide example.}
  \label{fig:slab-val-conv}
\end{figure}

As a concrete measure of the impact of the element size and polynomial order on the
value of $\frac{\lVert\bm{E}_{\mathrm{ref}}-\bm{E}_{\mathrm{PML}}\rVert}{\lVert\bm{E}_{\mathrm{ref}}\rVert}$,  Figure \ref{fig:slab-val-conv} presents the computed values for a range of
element sizes and varying polynomial order. It is clear that the obtained
value is depends heavily on the mesh and must be computed for each new analysis.

\subsection{On the needed number of modes}

The scenarios in which a WPBC is needed often involve reflections for which there is no analytical expression,
and for large problems, computing all the algebraic modes for the waveguide quickly becomes unpractical.
The last result of the previous subsection can be used as a tool for estimating if enough modes are being
employed in the WPBC: it is well known that Equation \ref{eq:slab-val-1}
is not expected to be satisfied at machine-level precision for unstructured meshes,
thus Equation \ref{eq:slab-val-res} can be seen as a tolerance level for analysing the efficiency of the WPBC.

In order to demonstrate it, a slab waveguide discontinuity is analysed, with core widths $w_1=\qty{0.4}{\micro\metre}$
and $w_2=\qty{1.5}{\micro\metre}$. On the slab with width $w_1$, the dominant mode is injected, and WPBCs are placed
on both sides \qty{1}{\micro\metre} away from the junction.
For this experiment, the number of modes used in each WPBC, $n_m^{\mathrm{in}}$ and $n_m^{\mathrm{out}}$, was varied from 3 to 300.

For each combination of $(n_m^{\mathrm{in}}$,$n_m^{\mathrm{out}})$, the following procedure was executed after
the scattering analysis was performed: at $\Gamma_{\mathrm{in}}^e$ (\emph{resp.} $\Gamma_{\mathrm{in}}^e$, now simply denoted by $\Gamma^e$),
a Finite Element space $\mathcal{U}$ was created and,
since the mesh nodes are periodic to $\Gamma_{\mathrm{in}}$ (\emph{resp.} $\Gamma_{\mathrm{out}}$), the computed
waveguide modes at the WPBC were transferred to $\mathcal{U}$.
Then, the 1D approximation space at $\Gamma^e$ was restricted such that its only free \emph{dofs} are the waveguide modes.
At the restricted $\mathcal{U}$, one can compute $u_p$, a $L^2(\Gamma^e)$-projection of the scattering solution $u$,
and the norm $r = \lVert u_p-u\rVert/\lVert u \rVert$ provides a measure on how well the employed modes are able to represent
the solution at a close distance to the WPBC. A big $r$ would indicate that there are higher-order modes being reflected
from the WPBC, and therefore more modes need be used at the boundary condition.

\begin{figure}[ht]
  \centering
  \includegraphics[width=0.8\linewidth]{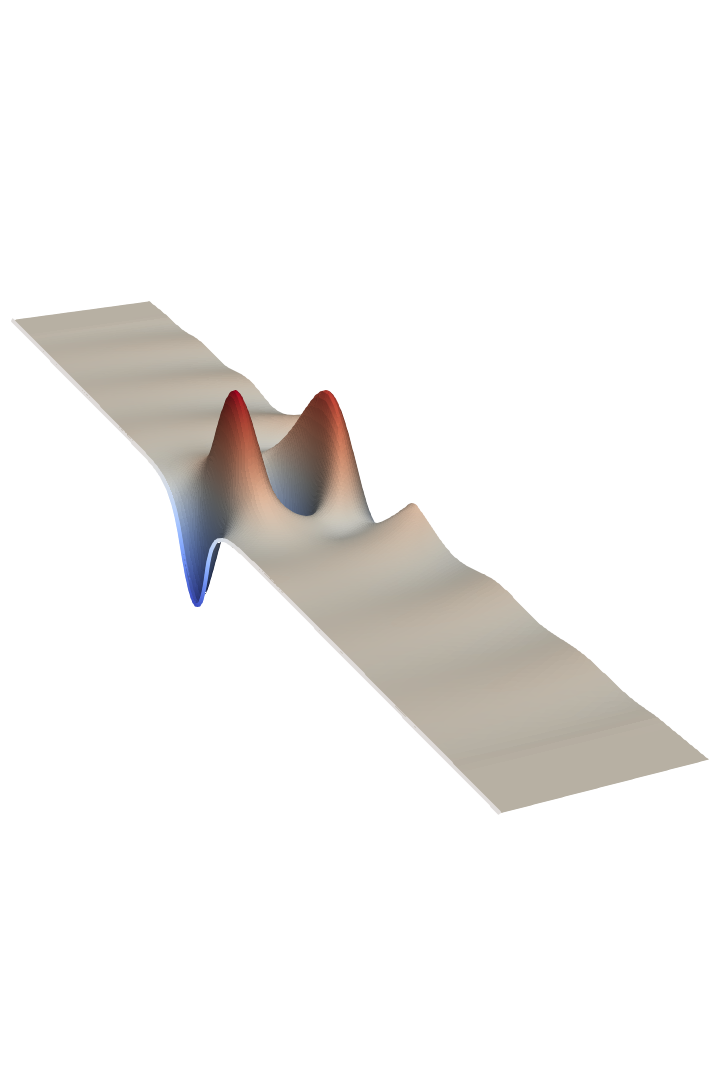}
  \caption{Dominant mode propagation in a slab waveguide discontinuity. Results obtained using the WPBC boundary condition on both domain truncations. The incident field is represented through a thicker line.}
  \label{fig:slab-disc-plot}
\end{figure}

\begin{figure*}[ht]
  \centering
  \includegraphics[width=0.8\linewidth]{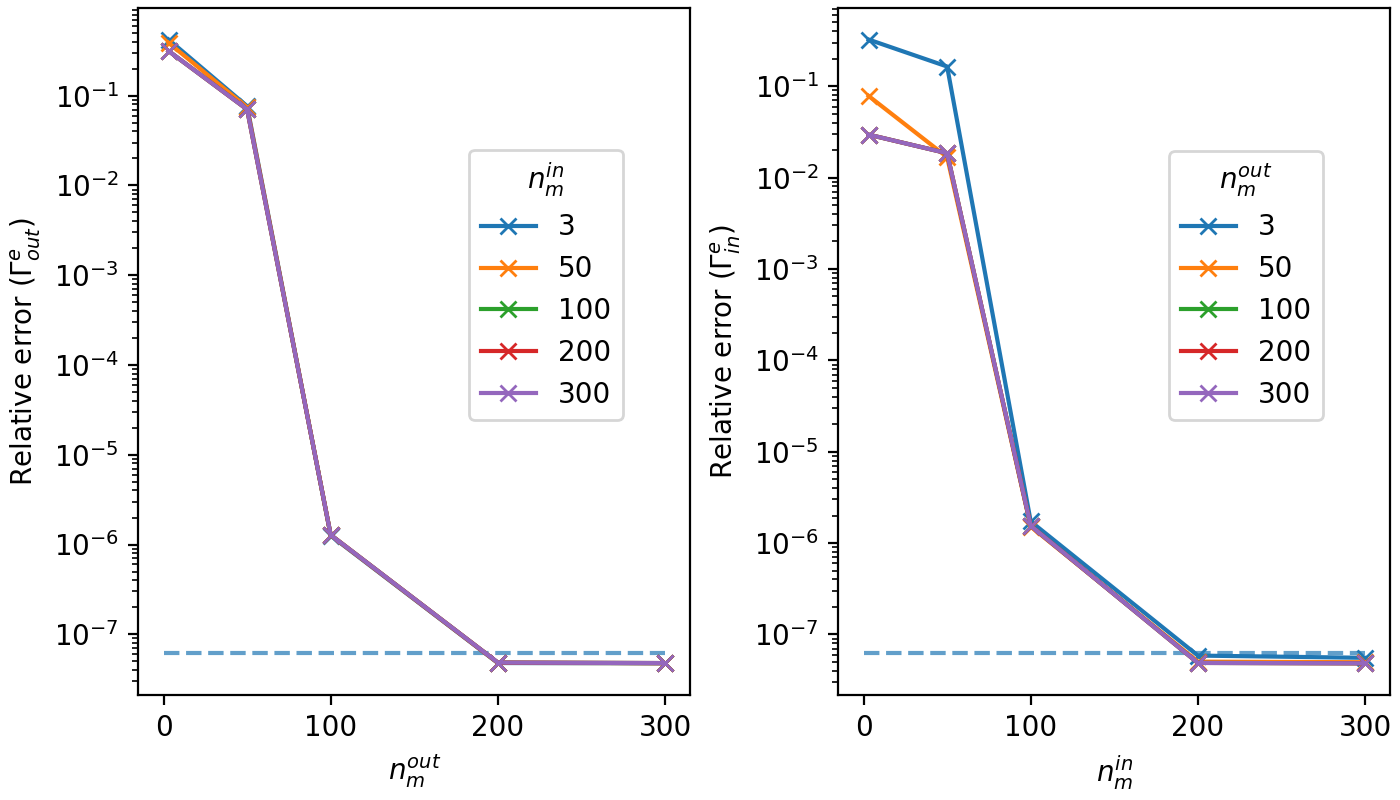}
  \caption{Relative error of the $\ltwo{\Gamma}$ projection of the scattering solution
    onto the waveguide modes at the evaluation subdomains with variable number of modes
    used at each WPBC for a slab waveguide discontinuity with
    slab widths $w_1=\qty{0.4}{\micro\metre}$ and $w_2=\qty{1.5}{\micro\metre}$}
  \label{fig:slab-disc-error}
\end{figure*}

Figure \ref{fig:slab-disc-error} shows the obtained $r$ for the slab discontinuity scattering shown in Figure \ref{fig:slab-disc-plot},
ran using the script \texttt{py/scripts/slab\_disc\_nmodes.py}.
The dashed line represents the numerical tolerance obtained from Equation \ref{eq:slab-val-res}: for this
analysis, 200 modes was considered sufficient on both domain truncations.

\subsection{Step fibre discontinuity}

In order to evaluate the competitive of the scheme in terms of computational time,
propagation along a step-index optical fibre discontinuity was used.
Due to the size of the linear system resulting from the 3D scattering analysis, a iterative solver
was used, namely GMRES, with a preconditioning scheme to be discuss in this section.

\begin{figure}[ht]
  \centering
  \includegraphics[width=0.8\linewidth]{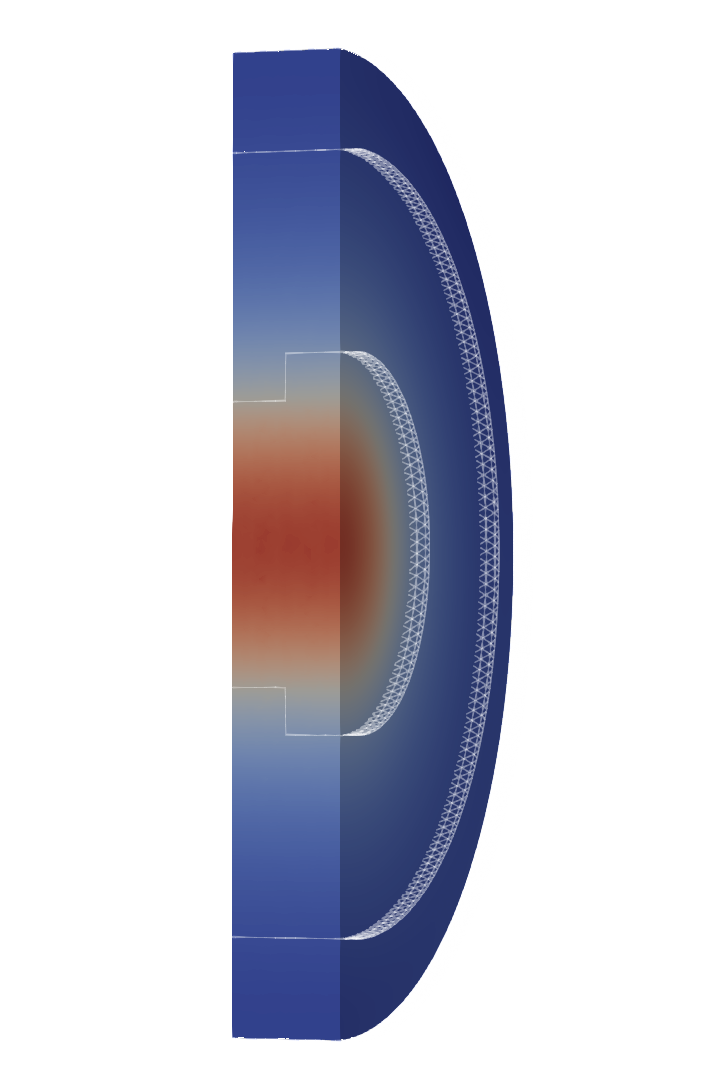}
  \caption{Geometry of the step-index optical fibre discontinuity problem superposed with
    the scattering solution obtained with the WPBC.}
  \label{fig:step-fibre-geo}
\end{figure}

The geometry of the problem is shown in Figure \ref{fig:step-fibre-geo}.
Element size was chosen as four elements per wavelength in the cladding and the
polynomial order as $k=3$.
This choice is due to limits on available computational resources, as the
memory requirements are larger when analysing 3D problems.
The fibres have $n_{\text{core}}=1.4457$ and $n_{\text{clad}}=1.4378$
and were analysed at $\lambda_0=\qty{1.55}{\micro\metre}$,
with the radial PML region placed at $3.5$ cladding wavelengths from the core.
At the chosen operational wavelength, even the dominant modes are not too confined
to the core, therefore the PML parameters were set at $\alpha_r = 0.4$ and, for the
PML-backed ports, $\alpha_z=1$, as to improve the conditioning of the matrix.

The coefficient
$\lVert\bm{E}_{\mathrm{ref}}-\bm{E}_{\mathrm{WPBC}}\rVert/\lVert\bm{E}_{\mathrm{ref}}\rVert$
was computed for a fibre with $r=\qty{6}{\micro\metre}$ and resulted in 0.00277218,
which, while considerably larger than the one in Equation \ref{eq:slab-val-res},
seems reasonable in the light of Figure \ref{fig:slab-val-conv}.
This computation was executed with the script \texttt{py/scripts/slab\_disc\_validation.py}.

\begin{figure*}[ht]
  \centering
  \includegraphics[width=0.8\linewidth]{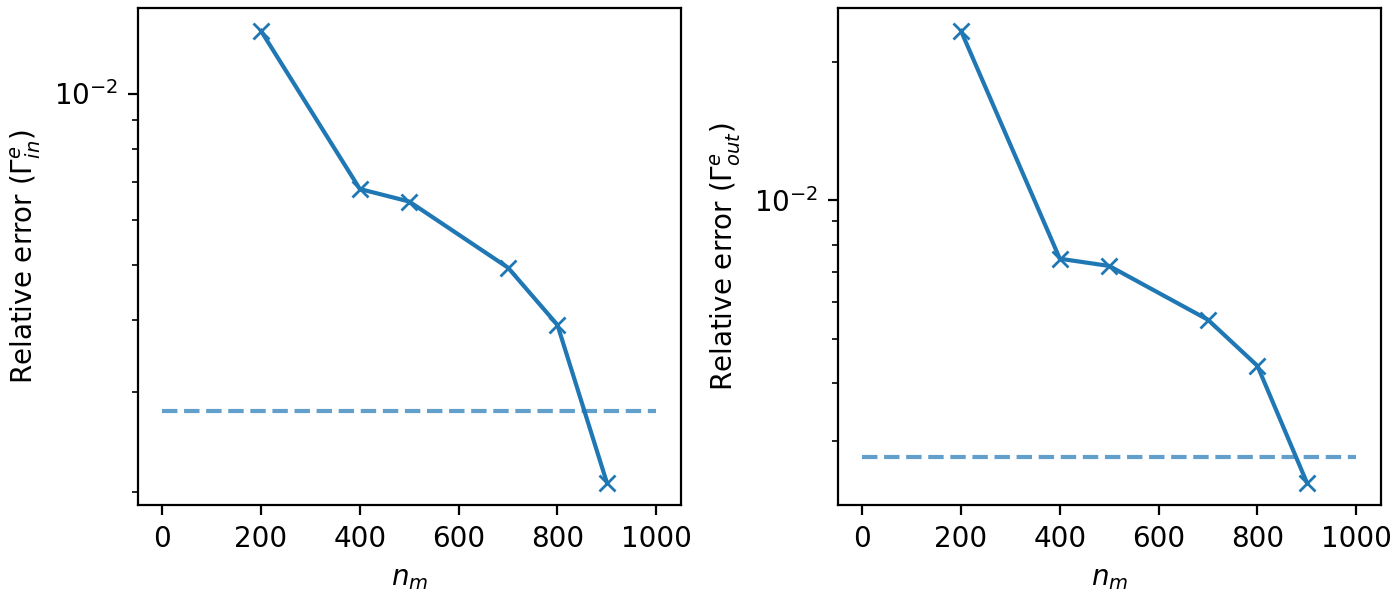}
  \caption{Relative error of the $\ltwo{\Gamma}$ projection of the scattering solution onto the waveguide modes at the evaluation subdomains with variable number of modes used at each WPBC for a step-index optical fibres discontinuity. The fibres had same refractive indices, and radii $r_1=\qty{6}{\micro\metre}$ and $r_2=\qty{8}{\micro\metre}$
    }
  \label{fig:sf3d-error}
\end{figure*}

The step-index discontinuity was setup with two fibres of radii
$r_1=\qty{6}{\micro\metre}$ and $r_2=\qty{8}{\micro\metre}$, and the waveguide ports
were placed at a distance of $2$ cladding wavelengths from the discontinuity.
Figure \ref{fig:sf3d-error} shows that 900 modes on each port were sufficient to
absorb the reflections due to the discontinuity,
therefore that was the chosen number of modes used in the comparison with
the PML-backed ports. 

For solving the resulting algebraic system of the scattering problem, a preconditioned
GMRES solver was chosen.
The preconditioning scheme was a block Jacobi preconditioner, with the following choice
of blocks, based on \citet{zaglmayr2006}:

\begin{equation}
  \{\mathcal{V}\}=\underbrace{\{\mathcal{E}_{l.o.}\}\oplus\{\mathcal{BC}\}}_{\text{sparse}}\oplus\{\mathcal{E}_{h.o.}\}\oplus\{\mathcal{F}\}\oplus\underbrace{\{\mathcal{I}\}}_{\text{c}}\text{.}
\end{equation}

The block $\{\mathcal{E}_{l.o.}\}$ is the block of the lowest order edge functions,
which is combined with $\{\mathcal{BC}\}$, the equations from the WPBC boundary,
as a sparse block to be decomposed by a direct solver. Then, the remaining blocks
are the small dense blocks $\{\mathcal{E}_{h.o.}\}$, associated with the higher order
edge functions, and $\{\mathcal{F}\}$, with the face functions.
The \emph{dofs} associated with the \emph{internal} functions of each element are
eliminated through static condensation, therefore they do not appear explicitly in the
algebraic system.
In Table \ref{table:sf3d-comp}, the results of the script \texttt{py/scripts/sf3d\_pml\_comp.py} are shown, comparing between the WPBC and PML-backed ports.
These were obtained in a computer with a Intel Xeon Gold 6130 CPU,
with 12 physical cores at \qty{2.1}{\giga\hertz} and \qty{376}{\giga\byte} of RAM memory.
The measured times are split into assembly time $t_{a}$, solving time $t_s$ and total time $t_t$.

One can see that, with reduced dimensions in the propagation direction and when a large number of modes is needed,
the assembly time of the WPBC setup is considerable.
However, the iterative solver converged in $40\%$ less iterations,
and the memory consumption was significantly reduced,
as the resulting computational domain was able to be shortened.

\begin{table*}[h]
    \begin{center}
      \begin{tabularx}{0.5\textwidth}{X Y Y}
        \toprule
                  &  PML     & WPBC      \\
        \midrule
        $n_{\textrm{dofs}}$ &  \num{26844552} & \num{3996635} \\
        $t_a(\unit{\second})$ &  \num{29319,518} & \num{6414,912}\\
        $t_s(\unit{\second})$ &  \num{6493,805} & \num{7291,170}\\
        $t_t(\unit{\second})$ &  \num{35813,323} & \num{13706,082}\\
        $n_{\textrm{iter}}$ &  \num{288} & \num{180}\\
      \bottomrule
    \end{tabularx}
  \end{center}
  \caption{Performance comparison of the scattering analysis of the step-index
    optical fibre discontinuity between the WPBC and the PML approaches. The times
  shown are the assembly, solving and total times.}
  \label{table:sf3d-comp}
\end{table*}

\subsection{Plasmonic sensor}

As a final demonstration of the capabilities of the WPBC,
the nanograting-based plasmonic sensor presented in \citep{he2022} is analysed.
The device is introduced in the context of biologic sensors, and presents a periodicity
that allows for the numeric analysis to be performed in a unit cell employing periodic
boundary conditions.

\begin{figure}[ht]
  \centering
  \includegraphics[width=\linewidth]{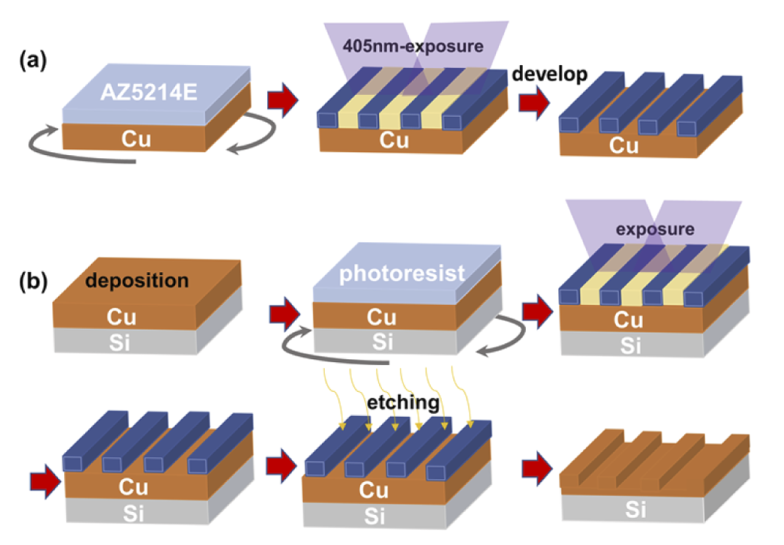}
  \caption{Geometry of the plasmonic sensor introduced in \citep{he2022}.
    For the numerical analysis of the unit cell, the left and right boundaries are set as Floquet periodic, controlling then the incidence angle of the light source. Image
  reproduced from \citep{he2022}.}
  \label{fig:plasmonic-geom}
\end{figure}

The geometry of the device is shown in Figure \ref{fig:plasmonic-geom}. The numerical
setup consists in illuminating the device from the top boundary. The left and right
boundaries are set to be Floquet periodic, in which a zero-valued phase coefficient
is associated with a normal incidence angle. In \citet{he2022}, the so-called DGC-SPR
sensor is compared with a simple copper grating.
Both devices consist of a copper layer with ridges. In the former case, the ridges
are made of the photoresist AZ5314E, in the latter, copper. The dielectric data is
presented in the supplementary material of the paper and the raw data was kindly
given to us upon request.

Even though the devices are not waveguides, the WPBC application is straightforward:
a one-dimensional modal analysis performed on the top boundary with periodic BCs
should provide the eigenfunctions for representing the reflected field produced by
illuminating the device, by the same concepts discussed in the introduction of Section
\ref{sec:wpbc}.

\begin{figure}[ht]
  \centering
  \includegraphics[width=\linewidth]{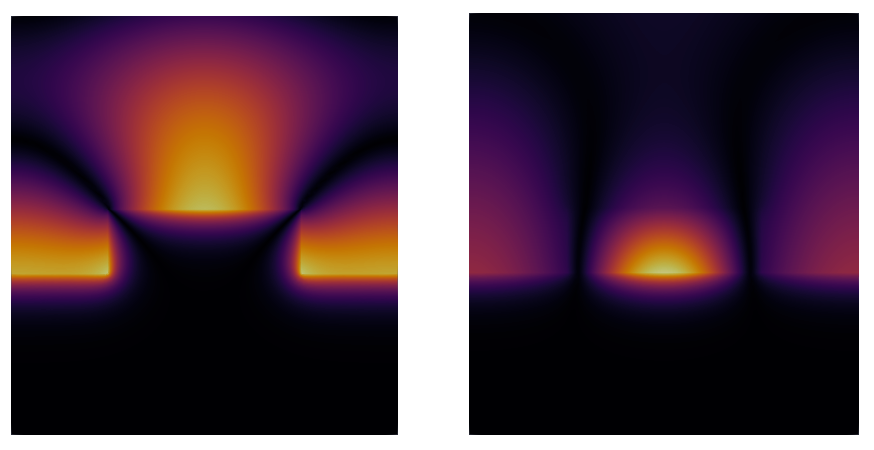}
  \caption{
    Field distributions obtained by illuminating the copper grating (left) and the
    DGC-SPR sensor(right) introduced in \citep{he2022} with wavelengths of
    $\lambda_0=\qty{746}{\nano\metre}$ and
    $\lambda_0=\qty{741}{\nano\metre}$, respectively.}
  \label{fig:plasmonic-field}
\end{figure}

In Figure \ref{fig:plasmonic-field}, the obtained field distributions corresponding
to Figures 2(e,f) of \citet{he2022} are shown, illustrating then the application
of the WPBC in such a scenario.

\section{Conclusions}

An implementation of the Waveguide Port Boundary Condition based on the restriction
of the approximation space is presented, avoiding the double integrals involved in the
original form of such boundary condition.
The implementation is validated in two and three dimensions along with the presentation
of a strategy for determining the needed number of modes for scenarios with strong
reflections, which are often the case in the numerical analysis of waveguide discontinuities.
A computationally intensive three dimensional example is used for a comparative
performance analysis with PML-backed waveguide ports, showing that despite the
considerably larger assembly time, it still presents advantage in terms of processing
time due to the reduction of the computational domain in the absence of PMLs.
The computation of the fields distribution of a nanograting-based plasmonic sensor is used
as an illustration of possible applications of the developed technique apart from
waveguide analysis. Future applications of the technique are under active development, with particular emphasis on its utilisation in the analysis of multiple sequential waveguide discontinuities.

\section{Acknowledgements}

Francisco T. Orlandini is grateful to the support during the development of this work
provided by the Brazilian Research Council (CNPq) (Grant 445074/2020-5).

The author Philippe R. B. Devloo is also grateful to CNPq (grants 305823/2017-5 and 306167/2017-4) and extends his gratitude to EPIC – Energy Production Innovation Center, hosted by the University of Campinas(UNICAMP) and sponsored by Equinor Brazil and FAPESP – São Paulo Research Foundation (2017/15736-3).

Hugo E. Hernández-Figueroa is thankful to the São Paulo Research Foundation (FAPESP) under the projects 2021/06506–0 (StReAM), 2021/11380–5 (CPTEn), 2021/00199–8 (SMARTNESS), and (EMU) 2022/11596-0; and the Brazilian Agency CNPq, under the project 314539/2023-9 (HEHF’s research productivity grant).

Finally, the authors of \citet{he2022} have kindly provided the dielectric data
for their plasmonic sensor allowing for us to analyse it as a final example in the
present work, and therefore deserve to be included here.

\begin{appendices}

\section{PML as an anisotropic absorber}\label{sec:pml}

The main idea of the Berenger PML is to introduce a complex coordinate stretching on a region
of the domain such that prescribed losses are slowly increased in a way that, for the
exact wave equation, incident waves are
absorbed in a reflectionless manner. The PML can also be seen as an anisotropic absorber\cite{sacks1995}, and
through the ideas in \citet{nicolet2008}, its equivalency can be shown.

Let us assume for now that the transformed coordinates inside the PML region are given by

\begin{equation}
  x' \rightarrow s_x x\text{,}\qquad y' \rightarrow s_y y\text{,}\qquad z' \rightarrow s_z z\text{,}
\end{equation}

where $s_i$, $i=x,y,z$ is the complex scaling coefficient causing the attenuation in the $i$ direction.
The coordinate transformation can be described by the Jacobian matrix

\begin{equation}
  \label{eq:pml-jac}
  \bm{J}_{\textrm{PML}} =
  \begin{pmatrix}
    s_x & 0 & 0 \\
    0 & s_y & 0 \\
    0 & 0 & s_z
  \end{pmatrix}\text{,}
\end{equation}
and the relative electric permittivity (magnetic permeability) $\bm{\epsilon}_r$ ($\bm{\mu}_r$)
can be transformed by \citep{nicolet2008}:

\begin{equation}
  \label{eq:pml-mat-transf}
  \bm{\gamma}_{\textrm{PML}} = \frac{1}{\det \bm{J}_{\textrm{PML}}}\bm{J}_{\textrm{PML}}^{-1}\bm{\gamma}\bm{J}^{-T}_{\textrm{PML}}\text{,}\quad \bm{\gamma} = \bm{\epsilon}_r, \bm{\mu}_r\text{,}
\end{equation}

which, for isotropic materials, result in

\begin{equation}
  \label{eq:pml-mat-transf-iso}
  \bm{\gamma}_{\textrm{PML}} = \gamma \bm{T}\text{,}\quad \bm{T} =
  \begin{pmatrix}
    \frac{s_y s_z}{s_x} & 0 & 0\\
    0 & \frac{s_x s_y}{s_y} & 0\\
    0 & 0 & \frac{s_x s_y}{s_z}
  \end{pmatrix}\text{.}
\end{equation}

Finally, for a cylindrical PML, the procedure is carried out similarly by composing the Jacobians
of the transformations to and from the cylindrical coordinate system.
\subsection{Determining the PML coefficients}\label{sec:pml-param}

In \citet{chiang2011}, the stretching coefficient $s_i$, $i=x,y,z$ is computed as
$s_i = \alpha_{\text{max}}^i (\frac{\rho_i}{d_i})^m$, where $\rho_i$ is the distance
from the beginning of the PML region, $d_i$ is the width of the PML region and

\begin{equation}
  \label{eq:pml-param}
  \alpha_{\text{max}}^i = -\frac{(m+1)\lambda}{4\pi n d_i}\ln{R}\text{,}
\end{equation}
with $n$ being the refractive index of the adjacent region and $R$ the desired
reflection coefficient. In the present work, $m=2$ and $R=10^{-70}$ were used
as initial parameters for setting up the PMLs.
\end{appendices}

\bibliography{references}

\end{document}